\documentclass[aps,prl,reprint,letterpaper,superscriptaddress,amsmath,showpacs,floatfix,longbibliography]{revtex4-2}
\RequirePackage{graphicx}
\usepackage{placeins}
\usepackage{hyperref}
\usepackage{xcolor}
\usepackage{amsfonts}

\renewcommand*{\div}{\ensuremath{\mathrm{div\,}}}
\newcommand*{\curl}{\ensuremath{\mathrm{curl\,}}}

\def\p{\partial}

\def\f1r{{\frac{1}{r}}  }

\begin{document}

\title{Asymptotic self-similar blow-up profile for three-dimensional
    axisymmetric Euler equations using neural networks}

\date{\today}

\author{Y. Wang}
\affiliation{Department of Geosciences, 
Princeton University
Princeton, NJ 08540, USA}

\author{C.-Y. Lai}
\affiliation{Department of Geosciences, 
Princeton University
Princeton, NJ 08540, USA}

\author{J. G\'omez-Serrano}
\affiliation{Department of Mathematics,
Brown University, Kassar House, 151 Thayer St., 
Providence, RI 02912, USA}
\affiliation{Departament de Matem\`atiques i Inform\`atica,
Universitat de Barcelona, Gran Via de les Corts Catalanes, 585,
08007, Barcelona, Spain}
\affiliation{
Centre de Recerca Matem\`atica,
Edifici C,
Campus Bellaterra,
08193 Bellaterra, Spain}

\author{T.~Buckmaster}
\email{Author to whom correspondence should be addressed.
Email: tbuckmaster@math.ias.edu}
\affiliation{School of Mathematics, Institute for Advanced Study
Princeton, NJ}
\affiliation{Department of Mathematics, Princeton University, Princeton, NJ 08544, USA}



\begin{abstract}\noindent 
Whether there exist finite time blow-up solutions for the 2-D Boussinesq and the 3-D Euler equations are of fundamental importance to the field of  fluid mechanics. We develop a new numerical framework, employing physics-informed neural networks (PINNs),  that discover, {\it for the first time}, a smooth self-similar blow-up profile for both equations. The solution itself could form the basis of a future computer-assisted proof of blow-up for both equations. In addition, we demonstrate PINNs could be successfully applied to find {\it unstable} self-similar solutions to fluid equations by constructing the first example of an unstable self-similar solution to the C\'ordoba-C\'ordoba-Fontelos equation.   We show that our numerical framework is both robust and  adaptable to various other equations.

\end{abstract}

\maketitle
A celebrated  open question in fluids is whether or not from smooth initial data the 3-D Euler equations may develop finite-time singularities (the inviscid analogue of the  Navier-Stokes Millennium-prize problem). For non-smooth $C^{1,\alpha}$ initial data with $0<\alpha\ll 1$, finite-time self-similar blow-up was proven in the groundbreaking work of Elgindi \cite{Elgindi2021, elgindi2019stability}.  The question of finite-time blow-up from {\it smooth} initial data remains {\it unresolved}. 

In the presence of a cylindrical boundary, Luo and Hou \cite{Luo12968} (cf.\ \cite{MR3278833}) performed compelling numerical simulations proposing a scenario -- sharing similarities with Pumir and Siggia \cite{PhysRevLett.68.1511} (cf.\ \cite{Childress,doi:10.1063/1.858422}) -- for finite time blow-up of the axi-symmetric 3-D Euler equations. They simulated the time-dependent problem and observed a dramatic growth in the maximum of vorticity (by a factor of $3 \cdot 10^{8}$), strongly suggesting formation of a singularity. The work is also suggestive of asymptotic self-similarity.

To confirm the existence of the finite-time singularity in the Luo-Hou scenario, and find its self-similar structure, we need to solve the self-similar equations associated with the axisymmetric 3-D Euler equations in the local coordinates near the singularity, which poses an extreme challenge to classical numerical methods as explained later. In this Letter, we develop a new numerical strategy, based on Physics-informed Neural Networks (PINNs)  that can solve the self-similar equations in a  simple and robust way. This new method allows us, for the first time, to find the smooth asymptotic self-similar blow-up profile for the Luo-Hou scenario. To the best of our knowledge, the solution is also the first truly multi-dimensional smooth backwards self-similar profile  for an equation from fluid mechanics.

Singularity formation for 3-D Euler equations with a cylindrical boundary is intrinsically linked to the same problem for the 2-D Bousinessq equations (cf.~\cite{MR3278833,Elgindi2020,ElJe2019,Chen2021}), another fundamental question in fluid mechanics, mentioned in Yudovich's \emph{`Eleven great problems of mathematical hydrodynamics'  }\cite{Yu03}. The mechanism for blow-up for the two equations is believed to be identical. The 2-D Boussinesq equations take the form
\begin{equation}\label{eq:Boussinesq}
\begin{gathered}
\p_t {\bf u}+ {\bf u}\cdot\nabla  {\bf u}+\nabla p= (0,\theta),\\\p_t \theta+  {\bf u}\cdot\nabla \theta=0,~\div  {\bf u}=0\,,
\end{gathered}
\end{equation}
where the 2-D vector $ {\bf u}( {\bf x},t)$ is the velocity and the scalar $\theta( {\bf x},t)$ is the temperature. We consider the spatial variable ${\bf x}=(x_1,x_2)$ to be taken on the half plane $x_2\geq 0$ and impose a non-penetration boundary condition at $x_2 = 0$ ($x_1$-axis), namely $u_2(x_1,0)=0$.

To search for singularity formation for the Boussinesq equations \eqref{eq:Boussinesq}, we look for \emph{backwards} \footnote{{Here \emph{backwards} self-similar solution refers to a solution that \emph{forms} a self-similar solution in finite time. Conversely, a \emph{forward} self-similar solution refers to a solution that starts from singular initial data and exists for all later time. This latter type of solution is for example relevant for studying the non-uniqueness of Leray-Hopf solutions to the Navier-Stokes equations.}} self-similar solutions of the form 
$
 {\bf u}=(1-t)^{\lambda} {\bf U}({\bf y})$ and $
\theta=(1-t)^{-1+\lambda}\Theta({\bf y})$, where we define the self-similar coordinates as
$
    {\bf y}=(y_1, y_2) = \tfrac{ (x_1, x_2)}{(1-t)^{1+\lambda}}
$
with ${\lambda>-1}$  yet to be determined. 
Under such a self-similar ansatz, the equations \eqref{eq:Boussinesq} become
 \begin{equation}\label{eq:SS:Boussinesq}
 \begin{gathered}
-\lambda  {\bf U}+((1+\lambda){\bf y}+ {\bf U})\cdot\nabla{\bf U}+\nabla P=(0,\Theta) \,,\\
(1-\lambda)\Theta +((1+\lambda){\bf y}+ {\bf U})\cdot\nabla \Theta=0,~\div {\bf U}=0 \,.
\end{gathered}
 \end{equation}
The corresponding solution is expected to have infinite energy \cite{Chae07}; however, we  impose that the solution to \eqref{eq:SS:Boussinesq} has mild growth at infinity (see cost functions in Supplementary Material) \footnote{{Specifically, we impose that the gradient of ${\bf U}$ vanishes at infinity, see cost functions in Supplementary Material.}}  which is an essential requirement for such a solution to be \emph{cut-off} to produce an asymptotically self-similar solution with a finite energy. \footnote{{Ideally, blow-up should be proven in the context of finite-energy solutions. To make use of an infinite-energy self-similar solution, one smoothly truncates the self-similar solution and then proves a form of stability. This procedure is standard. The enforced decay prohibits spurious solutions such as $({\bf U},\Theta)=(y,0)$.}}

Setting  $\Omega = \curl{\bf U}= \p_{y_1}U_2 -\p_{y_2}U_1$,  $\Phi=\p_{y_1}\Theta$ and $\Psi=\p_{y_2}\Theta$, we  rewrite \eqref{eq:SS:Boussinesq} in vorticity form
 \begin{equation}\label{eq:Kookaburra}
 \begin{split}
 \Omega+((1+\lambda){\bf y}+ {\bf U})\cdot\nabla \Omega&=\Phi \,, \\
(2+\partial_{y_1}U_{1})\Phi+((1+\lambda){\bf y}+  {\bf U})\cdot\nabla \Phi&=-\partial_{y_1} U_{2}\Psi \,, \\
 (2+\partial_{y_2}U_{2})\Psi+((1+\lambda){\bf y}+  {\bf U})\cdot\nabla \Psi&=-\partial_{y_2} U_{1}\Phi \,,\\
 \div {\bf U}&=0\, .
 \end{split} \end{equation}
To help find the solution, we impose the symmetries:  ($U_1, \Phi, \Omega$) are odd and $(U_2, \Psi)$ are even in the $y_1$ direction. In addition, we impose $ U_2(y_1,0)=0$ (the non-penetration boundary condition). To guarantee the uniqueness of the solution (removing scaling symmetry), we constrain
$\p_{y_1}\Omega (0, 0)=-1$. Finally, to rule out extraneous solutions, we impose that $\nabla  {\bf U}$, $\Phi$ and $\Psi$ all vanish at infinity. 

To describe the Luo-Hou scenario for the 3-D Euler blow-up in the presence of boundary, we  write the axisymmetric 3-D Euler equations: 
\begin{equation}\label{eq:Euler}
	 \begin{split}
		\left(\partial_t + u_r \partial_r+u_3 \partial_{x_3}\right)\left(\tfrac{\omega_\theta}{r}\right)
		&=\tfrac{1}{r^4} \partial_{x_3}\left(r u_\theta \right)^2\\
		\left(\partial_t + u_r \partial_r+u_3 \partial_{x_3}\right)\left(r u_\theta \right)&=0
	\end{split}
\end{equation}
where $(u_{r},u_\theta,u_{3})$ is the velocity in cylindrical coordinates and $\omega_\theta$ is the angular component of the vorticity. We  introduce a cylindrical boundary at $r=1$, and restrict  to the exterior domain \mbox{$\{(r,x_3) \in r \geq 1, x_3 \in \mathbb{R}\}$}. By imposing the self-similar ansatz
$(u_{x_3},u_r)=(1-t)^{\lambda} {\bf U}({\bf y}, s)$, $\omega_{\theta}=(1-t)^{-1}\Omega({\bf y},s)$, $
\p_r\left(r u_\theta \right)^2 =(1-t)^{-2}\Psi({\bf y},s)$ and $
\p_{x_3}\left(r u_\theta \right)^2=(1-t)^{-2}\Phi({\bf y},s)$
for self-similar coordinates 
${\bf y} =  (y_1,y_2)=\tfrac{(x_3,r-1)}{(1-t)^{1+\lambda}}$ and $s=-\log(1-t)$, the 3-D Euler equation \eqref{eq:Euler} becomes 
\begin{align}
({\partial_s}+\Omega)+((1+\lambda){\bf y} + {\bf U} )\cdot\nabla \Omega&=\Phi +\mathcal E_1\notag\\
({\partial_s}+2+\partial_{y_1}U_{1})\Phi+((1+\lambda){\bf y} +  {\bf U})\cdot\nabla \Phi&=-\partial_{y_1} U_{2}\Psi\notag \\
 ({\partial_s}+2+\partial_{y_2}U_{2})\Psi+((1+\lambda){\bf y} +  {\bf U})\cdot\nabla \Psi&=-\partial_{y_2} U_{1}\Phi \notag\\
 \div  {\bf U}& = \mathcal E_2\label{eq:SS:Euler}
\end{align}
where the errors $\mathcal E_1$ and $\mathcal E_2$ are given by the expressions
\begin{align*}\small
 \mathcal E_1&{\tiny= -y_2e^{-(1+\lambda)s}\tfrac{(y_2e^{-(1+\lambda)s}+2)(y_2^2e^{-2(1+\lambda)s}+2y_2e^{-(1+\lambda)s}+2)}{(1+y_2 e^{-(1+\lambda)s})^4} \Phi}\\
 \mathcal E_2 &{\tiny=-e^{-(1+\lambda) s}\tfrac{U_2}{1+y_2 e^{-(1+\lambda)s}}\,.}
\end{align*}
We look for solutions which are asymptotically self-similar: i.e.\ in self-similar coordinates they converge to a stationary state as $s\rightarrow \infty$. For such solutions, at any fixed ${\bf y}$, we have $\mathcal E_1,\mathcal E_2=O(e^{-(1+\lambda)s})$ and thus the errors decay exponentially fast in self-similar coordinates assuming that $\lambda>-1$. Thus, the self-similar equations \eqref{eq:SS:Euler} for Euler   converge to Bousinessq \eqref{eq:Kookaburra} as $s\rightarrow \infty$. Namely, the self-similar  solution for Boussinesq is identical to the asymptotic self-similar blow-up profile to Euler with cylindrical boundary.  

A key difficulty of solving the equation \eqref{eq:Kookaburra} lies in the unknown parameter $\lambda$ that needs to be solved simultaneously. We search for \emph{smooth} non-trivial solutions to \eqref{eq:Kookaburra}, which  exist for discrete $\lambda$ values.  This problem is extremely challenging using classical evolution (time-dependent) based numerical methods (cf.\ \cite{Cordoba-Fontelos-Mancho-Rodrigo:evidence-singularities-contour-dynamics,Scott-Dritschel:self-similar-sqg,Scott-Dritschel:self-similar-sqg-patch,Luo12968,MR3278833,Liu-Pego:singularities-water-waves,Eggers-Fontelos:self-similar-survey}).

\begin{figure}[t]
	\centering
	\includegraphics[width=3in]{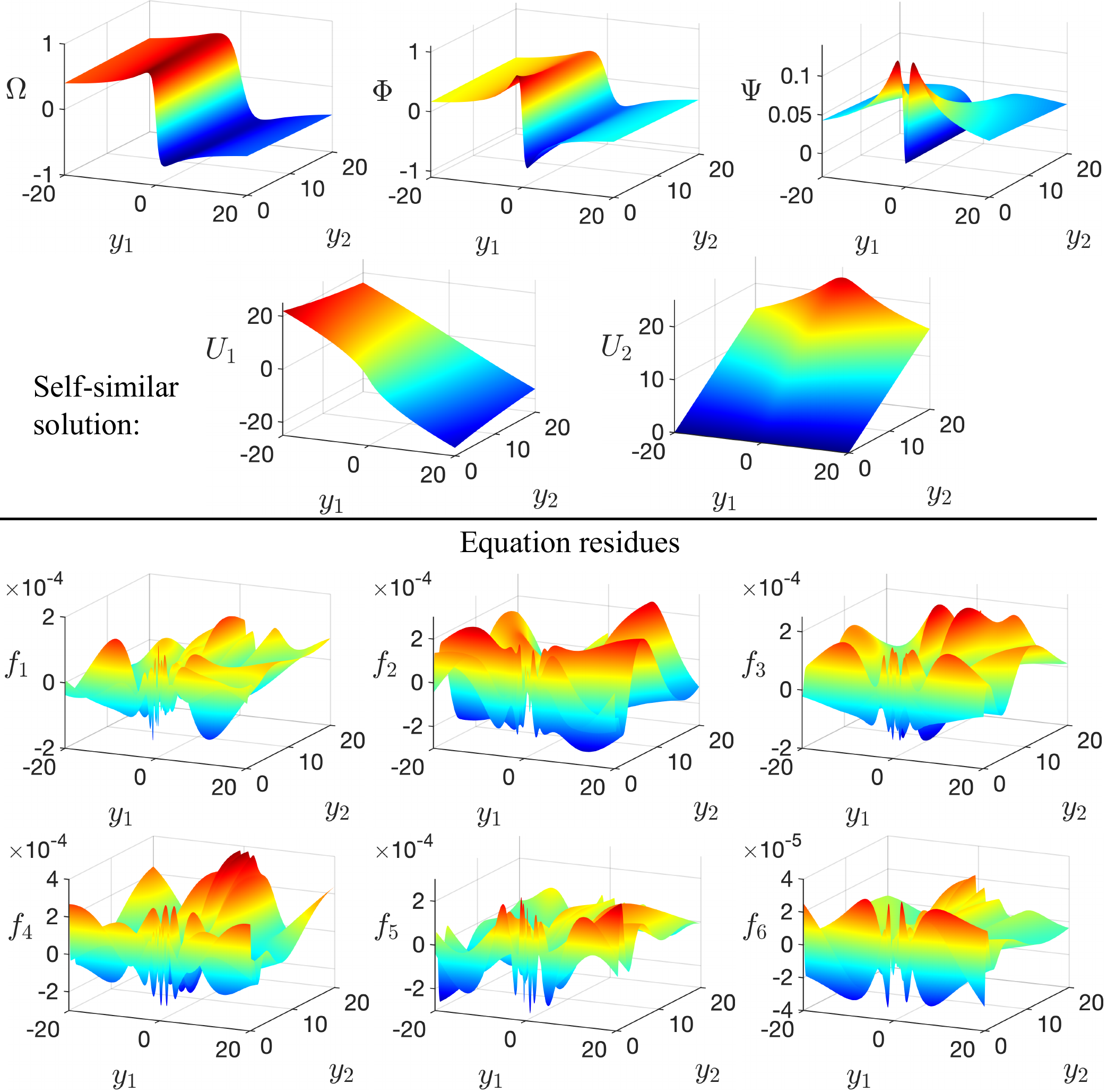}
	\caption{\label{fig:pdf} Solution for \eqref{eq:Kookaburra} found by PINN. $f_i$  indicate the residues, which are of five orders of magnitude smaller than  the solution. The inferred value of $\lambda$ is $1.917\pm 0.002$ after systematic test (see Supplementary Material).}
\end{figure}

Physics-informed neural networks (PINNs) were recently developed \ \cite{Raissi2019,Raissi2020} as a new class of numerical solver for PDEs and have been widely used in science and engineering \cite{Karniadakis2021}. In PINNs, neural networks approximate the solution to a PDE by searching for a solution in a continuous domain that approximately satisfies the physics constraints (e.g. equations and solution constraints).  PINNs have been successfully used to solve not only forward problems but also inverse problems (e.g.\ identifying the Reynolds number from a given flow and the Navier-Stokes equation \cite{Raissi2019}), demonstrating the capacity of PINNs to invert for unknown parameters in the governing equations. Here, we use a PINN to find not only the self similar solution profile but also the unknown self-similarity exponent $\lambda$. To guarantee the success of PINN, it is critical to understand the key symmetries of the problem, its spurious solutions, as well as intuition of the qualitative properties of the solution (e.g.\ its geometry and asymptotics).

To find the self-similar solution for the Bousinessq equation, we represent each of $U_1$, $U_2$, $\Omega$, $\Phi$ or $\Psi$ in the Bousinessq equations \eqref{eq:Kookaburra} by an individual fully-connected neural network with $y_1$ and $y_2$ as its inputs. We use 6 hidden layers with 30 units in each hidden layer for each network and use the hyperbolic tangent function $\tanh$ as the activation function.   We impose the symmetry of each variable ($U_1, U_2, \Omega, \Phi, \Psi$) by constructing the function form $q_{odd} = [\mathrm{NN}_q(y_1, y_2) - \mathrm{NN}_q(-y_1, y_2)]/2$ and $q_{even} = [\mathrm{NN}_q(y_1, y_2) + \mathrm{NN}_q(-y_1, y_2)]/2$, where $\mathrm{NN}_q$ is the neural network created for the variable $q$.   

To train the neural network, we need a cost function  and an optimization algorithm. For PINNs,  the cost function is composed of two types of loss. The first is the {\it condition loss}, which evaluates the residue of the solution condition, where the residue here is defined as the difference between the neural network approximated condition and the true solution condition. The condition loss can be written as
\begin{equation}\label{eq:loss_d}
	loss_c^{(j)} = \frac{1}{N_c^{(j)} }\sum_{i=1}^{N_c^{(j)} }g_{(j)}^2\left[ {\bf y}_i, \hat{q} \left({\bf y}_i \right)\right],
\end{equation}
where $g_{(j)} \left({\bf y}_i, \hat{q}({\bf y}_i) \right)$ indicates the residue of the $j$-th boundary condition at the $i$-th position ${\bf y}_i = (y_1, y_2)_i$ and $\hat{q}\left({\bf y}_i \right)$ indicates the neural network prediction of the variable $q$. The parameter $N_c^{(j)}$ indicates the total number of points used for evaluating the $j$-th boundary condition. 

The second type of loss is known as the {\it equation loss}, which evaluates the residue of the governing equation averaged over a set of collocation points over the domain.  The residue of the equation $f_{(k)}$ is defined as the error in the equation calculated with the neural network predictions. The equation loss can be written as
\begin{equation}\label{eq:loss_f}
loss_{f}^{(k)}  = \tfrac{1}{N_f^{(k)}}\sum_{i=1}^{N_f^{(k)}}f_{(k)}^2 \left[{\bf y}_i, \hat{q} \left({\bf y}_i \right) \right],
\end{equation}
where $f_{(k)}\left({\bf y_i}, \hat{q}({\bf y}_i) \right)$ indicates the residue of the $k$-th equation evaluated at the $i$-th collocation point.  The parameter $N_f^{(k)}$ denotes the total number of collocation points used for the $k$-th equation.   The residues of the boundary conditions $g_{(j)}$ and equations $f_{(k)}$ involved in the cost function are listed in the Supplemental Material.

\begin{figure}[t]
	\centering
	\includegraphics[width=3.3in]{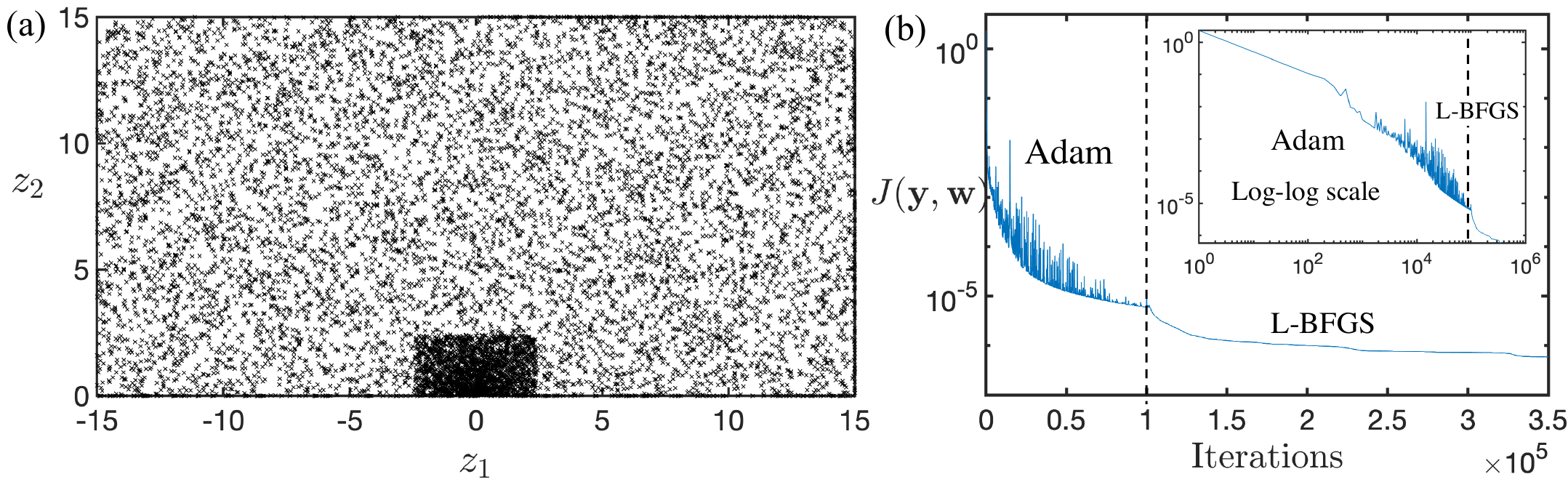}
	\caption{\label{fig:collo} (a) Spatial distribution of the collocation points. Here ($z_1$, $z_2$) are the rescaled coordinates: $y_1 = \sinh(z_1)$ and $y_2 = \sinh(z_2)$. 10,000 collocation points in total are used for training. (b) Decrease of the total loss over the training iterations.  The inset shows the loss curve in  a log-log scale.}
\end{figure}

We stress that all equations are local, which is an advantage of our method versus alternate methods that require a careful consideration of non-locality in infinite domains. In our implementation, to approximate an infinite domain, we introduce the coordinates $z=(z_1,z_2) = (\sinh^{-1}(y_1),\sinh^{-1}(y_2))$ (see Figure \ref{fig:collo}$a$) and consider a domain $z\in [-30,30]^2$. In the  $y$-coordinates, this corresponds to a domain $\approx [-5\cdot 10^{12},5\cdot 10^{12}]^2$. The equations written in $z$-coordinates are given in the Supplemental Material.

The equations and conditions provided so far can only find the unique solution to \eqref{eq:Kookaburra} for a specified $\lambda$. Figure 1 in the Supplementary Material shows the PINN solution to \eqref{eq:Kookaburra} for $\lambda = 3$. The large equation residue at the origin indicates the non-smoothness of the solution at the origin (see Supplemental Material). To search for the right $\lambda$ that guarantees the smoothness of solutions, i.e. avoiding the local peaks in the equation residue, we impose the \textit{smoothness constraint} to penalize the gradient of equation residues around the non-smooth point (origin). 
\begin{equation}\label{eq:loss_g}
loss_{s}^{(k)}  = \frac{1}{N_s^{(k)}}\sum_{i=1}^{N_s^{(k)}} \left|\nabla f_{(k)} \left({\bf y}_i, \bar{q} \left({\bf y}_i \right) \right) \right|^2\, ,
\end{equation}
where $N_s^{(k)}$ indicates the total number of collocation points close to the origin.  The sum of \eqref{eq:loss_d}, \eqref{eq:loss_f} and \eqref{eq:loss_g} defines the final cost function 
\[
	J({\bf y, w}) = \tfrac{1}{n_b}\sum_{j=1}^{n_b} loss_{c}^{(j)} + \gamma ( \tfrac{1}{n_e}\sum_{k=1}^{n_e} loss_{f}^{(k)} +\tfrac{1}{n_e}\sum_{k=1}^{n_e} loss_{s}^{(k)} )
\]
where $n_b=8$ and $n_e=6$ are the total number of solution conditions and governing equations used (see Supplementary Material).  The constant $\gamma$ is a hyper-parameter of PINNs, known as the {\it equation weight} \cite{van2021}, which balances the contribution of the condition loss and equation loss in the final cost function $J({\bf y, w})$. For Boussinesq, we choose $\gamma = 0.1$ for optimal training performance.

The common optimization methods used for PINN training are {\it Adam} \cite{kingma2014} and {\it L-BFGS} \cite{liu1989}. Despite the fact that no optimization method guarantees the convergence to a global minimum, our empirical experience, consistent with prior studies \cite{markidis2021},  shows that Adam performs better at avoiding local minima, while L-BFGS has a faster convergence rate throughout the training.  Thus, we use Adam first for 100,000 iterations and then L-BFGS for 250,000 iterations to search for the self-similar solution for the Bousinessq equations. Figure \ref{fig:collo}(b) shows the convergence of the cost function $J({\bf y, w})$ throughout the training iterations.

Figure \ref{fig:pdf} shows the approximate solutions to \eqref{eq:Kookaburra}, along with their corresponding equation residues. The equation residues are approximately five orders of magnitude smaller than that of the solution found. With the smoothness constraint \eqref{eq:loss_g} the inferred exponent for the smooth solution is $\lambda\approx1.917$. In the Supplementary Material we show the robustness of the PINN prediction with different random initialization and normalization condition. We also demonstrate convergence of the inferred $\lambda$ with domain size. The solutions found by the PINN are in agreement with the asymptotics of the time dependent solutions found by  Luo and Hou \cite{Luo12968,MR3278833}. Extrapolating from the paper \cite{MR3278833}, the work is suggestive of a self-similarity exponent of $\lambda \approx 1.9$, in agreement with the exponent found by the PINN. Similarly to Luo and Hou, the trajectories corresponding to the self-similar velocity follow the geometry of a hyperbolic point at the origin. 

The spatial distribution of the collocation points plays a critical role for the success of PINN training. To guide the neural network to find the correct self-similar solution for Bousinessq \eqref{eq:Kookaburra}, we train the neural network to prioritize the equation constraints around the origin. Towards  this goal, we divide the domain into two regions, one close to the origin and one far, in each region the collocation points are uniformly distributed. We increase the number of collocation points surrounding the origin. Otherwise, the neural network prediction would likely be trapped in a local minimum during the training.

The PINN-based scheme for finding self-similar blow-up offers advantages in terms of both universality and efficiency.  For universality, the above PINN scheme can be generally applied to solving various self-similar equations without the requirement of prior knowledge of specific structure. For efficiency, the smooth self-similar solution was in fact found by PINN throughout one single training. There is no continuation scheme or time evolution required for the training, largely reducing the computational cost of the method.  An additional major advantage of the PINN scheme, as we will later demonstrate, is its ability to find unstable self-similar solutions which would be incredibly difficult, if not impossible to find via traditional methods.

To validate our approach, we compare self-similar solutions obtained to known results in the literature (Supplementary Material). We apply the PINN scheme to find non-smooth solutions to the Boussinesq equation, which are in agreement with the explicit approximate solutions of Chen and Hou \cite{Chen2021} (Section 1.5 of the Supplementary Material).

One of the simplest PDEs  exhibiting self-similar blow-up is the 1-D Burgers' equation, which can be solved analytically. The equation provides an excellent sandbox to test and refine the PINN. Section 2 of the Supplementary Material shows that the PINN scheme can find stable, unstable and non-smooths self-similar solutions to the Burgers' equation. A common numerical strategy to finding self-similar solutions is to introduce time dependence into the problem: while this is straightforward in the stable case, instabilities in the unstable case make finding unstable self-similar profiles comparatively more difficult, if not impossible. The PINN  scheme does not suffer this drawback and thus presents itself as a great method for finding unstable smooth self-similar solutions. This later fact will be reinforced below where we demonstrate that the PINN is successful in finding new unstable self-similar solutions that have applicability to an important open problem in mathematical fluid dynamics.

\begin{figure}[t]
	\centering
	\includegraphics[width=3.4in]{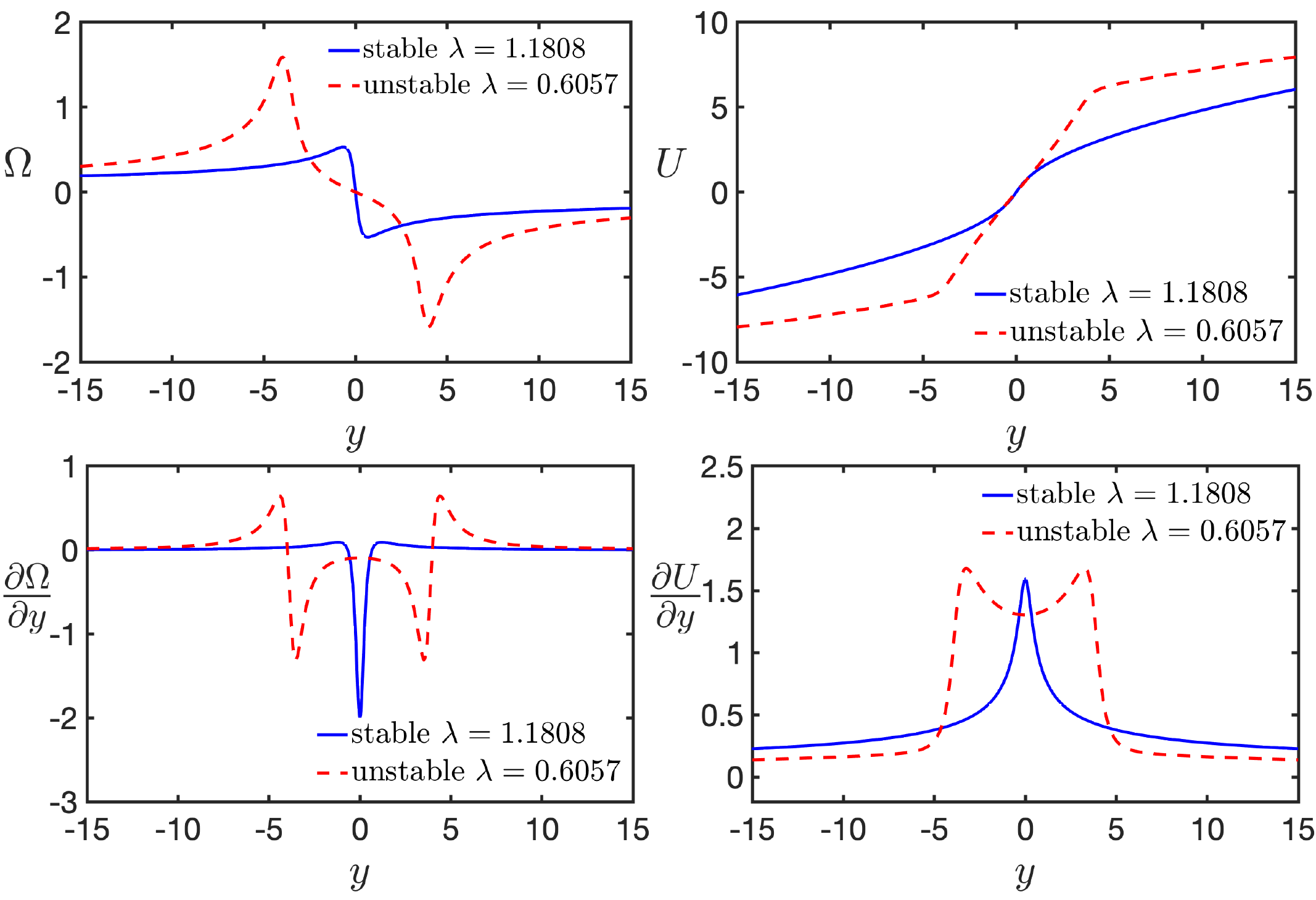}
	\caption{\label{fig:unstable} Comparison between stable and unstable solutions with associated $\lambda$ for CCF equation.  }
\end{figure}

The generalized De Gregorio equation  \cite{OkSaWu08} is given by
\begin{equation*}
\omega_t+au\omega_x=\omega u_x,\quad\mbox{where } u=\smallint_0^x H\omega=\Lambda^{-1}\omega
\end{equation*}
and $H$ is the Hilbert transform. The equation is a generalization of the De Gregorio equation  ($a=1$) \cite{DeGregorio1996}  and has been proposed as a one-dimensional model for an equation for which there is nontrivial interaction of  advection and vortex stretching (modeling behavior of the 3-D Euler equations). 

The case $a=0$, in the absence of advection,  is known in the literature as the Constantin-Lax-Majda equation. In this simple case, exact self-similar blow-up solutions can be constructed \cite{MR812343}. The case $a = -1$ (known as the C\'ordoba-C\'ordoba-Fontelos (CCF) model) was proposed  as a model of the surface quasi-geostrophic (SQG) equation and it also develops finite time singularities  \cite{CCF05}. In the case $a < 0$, advection and vortex stretching  work in conjunction leading to finite time   singularities \cite{CaCo10}. The case $a > 0$ leads to the competition of the two terms. By a clever expansion in $a$, smooth self-similar profiles were constructed in \cite{Elgindi2019} for small, positive $a$, leading to finite time blow-up. Via a computer-assisted proof, Chen, Hou, and Huang in \cite{ChHoHu19} proved blow-up for the De Gregorio equation  ($a=1$).

Lushnikov et al.\ \cite{Lushnikov2021} represents the most thorough numerical study of the generalized De Gregorio equation to date. In \cite{Lushnikov2021},  self-similar solutions were found in the whole range $a \in [-1,1]$ and  beyond. We use their reported parameters as benchmarks for our results. In the Supplementary Material \footnote{In Section 3 of the Supplementary material we reproduce the finds of \cite{Lushnikov2021} for the generalized De Gregorio equation.}, we show that PINN can accurately reproduce the findings of \cite{Lushnikov2021}. 

Returning to the specific case of CCF ($a=-1$), an interesting question to add fractional dissipation $(-\Delta)^{\alpha/2}$ and ask for which values of $\alpha$ do singularities occur. Blow-up is known to occur for $0
\leq \alpha <\frac12$, whereas for $\alpha\geq 1$ the problem is global-wellposed \cite{Li2008,MR2464570,Kiselev2010}. The behavior in the range $\frac12
\leq \alpha<1$ remains an important open problem.

Assuming the ansatz $\omega=\tfrac{1}{1-t}\Omega(\tfrac{x}{(1-t)^{1+\lambda}},s)$ for $s=-\log(1-t)$, then  in the self-similar evolution equation for $\Omega$, the dissipative term takes the form $e^{((1+\lambda)\alpha-1)s}(-\Delta)^{\frac\alpha2}\Omega$. Analogous to how self-similar Boussinesq solutions can be used to construct asymptotically self-similar solutions to Euler with boundary, self-similar inviscid CCF solutions satisfying the condition $ (1+\lambda)\alpha-1<0$ may be employed to construct asymptotically self-similar solutions to dissipative CCF. Since $\lambda\approx 1.18078$ for the stable self-similar solution to inviscid CCF, such a solution is  ill-suited to prove blow-up in the parameter range $\frac12 \leq \alpha<1$. Motivated by known work on the Burgers' equation \cite{EgFo2009,PhysRevFluids.3.110503} and compressible Euler \cite{MeRaRoSz19a, implosion}, one could conjecture the existence of a discrete hierarchy of unstable solutions with decreasing $\lambda$. By windowing the parameter $\lambda$, including additional derivatives of the governing equation in our residues, and using the constraint $\Omega(0.5)=0.05$ to renormalize, the PINN discovers  an unstable self-similar solution corresponding to $\lambda\approx 0.60573$ (see Figure \ref{fig:unstable}). Such a solution would allow us to prove blow-up for dissipative CCF for the range $\alpha<\frac{1}{1+\lambda}\approx 0.61$. Moreover, such a result is suggestive of a possible strategy of addressing the Navier-Stokes \emph{Millennium Prize} \cite{Fefferman00}, i.e.\ via unstable self-similar solutions to 3-D Euler (the same strategy has proved successful for dissipative Burgers' and compressible Navier-Stokes \cite{OhPa21,MeRaRoSz19c, implosion}). One  expects that the PINN may be adapted to find higher order unstable solutions to CCF as well as unstable solutions to the Boussinesq equation -- this is subject of future work.

\begin{acknowledgments}
The authors were supported by the NSF grants DMS-1900149 and DMS-1929284, ERC grant 852741, the Spanish State Research Agency, through the Severo Ochoa and Mar\'ia de Maeztu Program for Centers and Units of Excellence in R\&D (CEX2020-001084-M), Dean for Research Fund at Princeton University, the Schmidt DataX Fund at Princeton University,  a Simons Foundation Mathematical and Physical Sciences Collaborative Grant and grant from the Institute for Advanced Study. The simulations presented in this article were performed using the Princeton Research Computing resources at Princeton University and School of Natural Sciences Computing resources at the Institute for Advanced Study. 
\end{acknowledgments}

\nocite{bilato2014algorithm}

\end{document}


\maketitle
\vspace{-2em}
\section{The Boussinesq equations}
\subsection{Coordinate transformation}\label{sec:cord}
For the self-similar Bousinessq equations
\begin{equation}\label{eq:Kookaburra}
	\begin{split}
		\Omega+((1+\lambda){\bf y}+ {\bf U})\cdot\nabla \Omega&=\Phi \\
		(2+\partial_{y_1}U_{1})\Phi+((1+\lambda){\bf y}+  {\bf U})\cdot\nabla \Phi&=-\partial_{y_1} U_{2}\Psi \\
		(2+\partial_{y_2}U_{2})\Psi+((1+\lambda){\bf y}+  {\bf U})\cdot\nabla \Psi&=-\partial_{y_2} U_{1}\Phi \\
		\Omega = \partial_{y_2}U_{1} - \partial_{y_1}U_{2}  \qquad  \text{and} \qquad  \div & {\bf U}=0\ .
\end{split} \end{equation}

\begin{figure}[t]
	\centering
	\includegraphics[width=\textwidth]{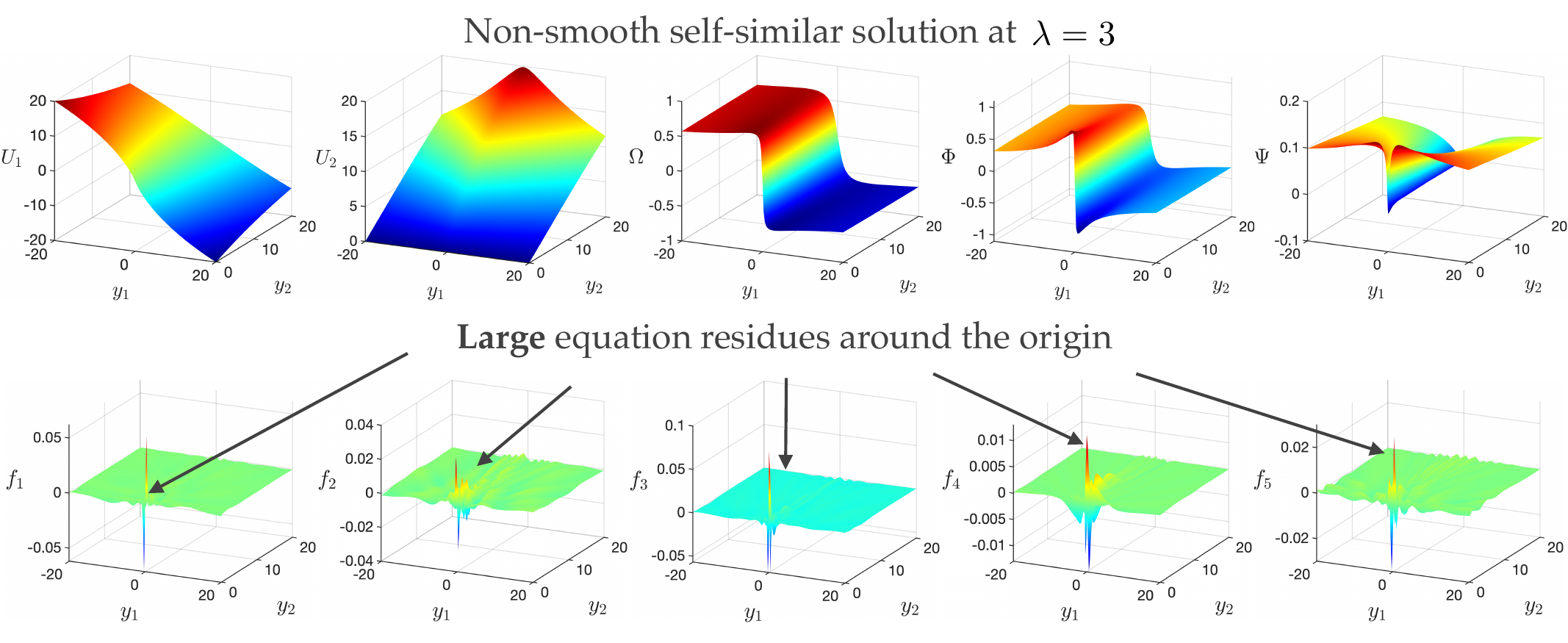} \vspace{-2em}
	\caption{\label{fig:bsqns} Non-smooth solutions to the 2D Boussinesq equations \eqref{eq:Kookaburra} derived by the physics-informed neural network for $\lambda = 3$, and its corresponding equation residues. }
\end{figure}

\noindent In order to cope with very large domains, we defined a new $z$-coordinate ($z_1$, $z_2$) with relation
\[{\bf y}=(y_1,y_2) = (\sinh(z_1), \: \sinh(z_2))  \quad \Longleftrightarrow  \quad {\bf z} =(z_1, z_2) = (\sinh^{-1} y_1, \: \sinh^{-1} y_2)   \,.\] 
which leads us to tackle domains with moderate dimensions in ${\bf z}$, but at the same time exponentially large in ${\bf y}$. The self-similar Boussinesq equations in ${\bf z}$ coordinates can be written as:
\begin{equation}\label{eq:Kookaburra_z}
	\begin{split}
		\Omega+((1+\lambda)( \sinh(z_1),\sinh(z_2))+ {\bf U})\cdot\nabla_y \Omega&=\Phi \\
		(2+\sech(z_1)\partial_{z_1}U_{1})\Phi+((1+\lambda)( \sinh(z_1),\sinh(z_2))+  {\bf U})\cdot\nabla_y \Phi&=-\sech(z_1)\partial_{z_1} U_{2}\Psi \\
		(2+\sech(z_2)\partial_{z_2}U_{2})\Psi+((1+\lambda)( \sinh(z_1),\sinh(z_2))+  {\bf U})\cdot\nabla_y \Psi&=-\sech(z_2)\partial_{z_2} U_{1}\Phi \\
		\Omega = \sech(z_2)\partial_{z_2}U_{1} - \sech(z_1)\partial_{z_1}U_{2}   \qquad  \text{and} \qquad  & \nabla_y \cdot {\bf U}  = 0
\end{split} \end{equation}
where $\nabla_y$ indicates the gradient operator with respect to the original $y$-coordinate, which can be expressed in the $z$-coordinate by $\nabla_y = (\sech(z_1) \partial_{z_1}, \: \sech(z_2) \partial_{z_2})$

\subsection{PINN cost functions}
The \textit{condition loss} and \textit{equation loss} included in the cost function are defined below. The \textit{condition loss}, which evaluates the residue of the boundary conditions includes
\begin{eqnarray*}
	g_1 &=&  U_2(y_1,0),   \qquad \ \: \text{ (non-penetration condition) }\\
	g_2 &=& \partial_{y_1}\Omega (0, 0) + 1,  \  \text{ (fix scaling symmetry) }\\
	g_3 &=& \nabla  {\bf U}(|y_1|=d),  \quad g_4 = \nabla  {\bf U}(|y_2|=d) \quad \text{(suitable decay)} \\
	g_5 &=& \Phi(|y_1|=d),  \quad g_6 = \Phi(|y_2|=d),\\
		g_7 &=& \Psi(|y_1|=d),  \quad g_8 = \Psi(|y_2|=d),
\end{eqnarray*}
where $d$ is the domain size of calculation.  The condition loss in the $z$-coordinates becomes
\begin{equation} \label{eq:bsqg}
	\begin{split}
	g_1 &=  U_2(z_1,0),   \qquad \ \: \text{ (non-penetration condition) }\\
	g_2 &= \partial_{z_1}\Omega (0, 0) + 1,  \  \text{ (fix scaling symmetry) }\\
	g_3 &= \nabla _y {\bf U}(|z_1|=\sinh^{-1} d),  \quad g_4 = \nabla_y  {\bf U}(|z_2|=\sinh^{-1}d) \quad \text{(suitable decay)} \\
	g_5 &= \Phi(|z_1|=\sinh^{-1}d),  \quad g_6 = \Phi(|z_2|=\sinh^{-1}d),\\
	g_7 &= \Psi(|z_1|=\sinh^{-1}d),  \quad g_8 = \Psi(|z_2|=\sinh^{-1}d),
	\end{split}
\end{equation}
where $\nabla_y = (\sech(z_1) \partial_{z_1}, \: \sech(z_2) \partial_{z_2})$. All governing equations for Bousinessq in the $z$-coordinates are listed in \eqref{eq:Kookaburra_z}. Thus, the \textit{equation loss} evaluates the residue of 6 equations in total:
\begin{equation} \label{eq:bsqf}
	\begin{split}
		f_1 &= \Omega+((1+\lambda)( \sinh(z_1),\sinh(z_2))+ {\bf U})\cdot\nabla_y \Omega - \Phi\\
		f_2 &= (2+\sech(z_1)\partial_{z_1}U_{1})\Phi+((1+\lambda)( \sinh(z_1),\sinh(z_2))+  {\bf U})\cdot\nabla_y \Phi + \sech(z_1)\partial_{z_1} U_{2}\Psi \\
		f_3 &= (2+\sech(z_2)\partial_{z_2}U_{2})\Psi+((1+\lambda)( \sinh(z_1),\sinh(z_2))+  {\bf U})\cdot\nabla_y \Psi + \sech(z_2)\partial_{z_2} U_{1}\Phi\\
		f_4 &= \sech(z_1)\partial_{z_1} U_{1} + \sech(z_2)\partial_{z_2} U_{2}, \\
		f_5 &= \Omega - \left[\sech(z_2)\partial_{z_2}U_{1} - \sech(z_1)\partial_{z_1}U_{2} \right], \\
		f_6 &=\sech(z_2)\partial_{z_2} \Phi  - \sech(z_1)\partial_{z_1} \Psi \: .
	\end{split}
\end{equation}

\begin{figure}[t]
	\centering
	\includegraphics[width=\textwidth]{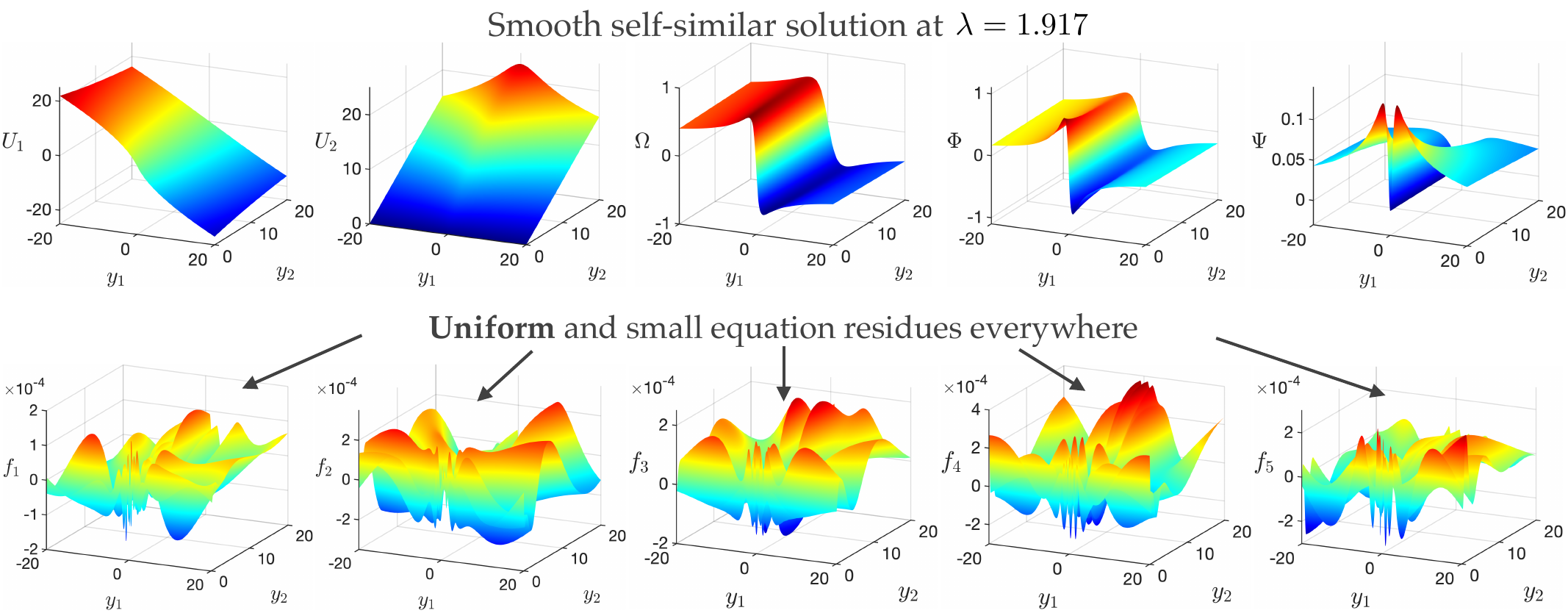}
	\caption{\label{fig:bsqsol} Smooth solutions to the 2D Boussinesq equations \eqref{eq:Kookaburra} derived by the physics-informed neural network with inferred $\lambda = 1.917$, and its corresponding equation residues. }
\end{figure}

\subsection{Smoothness constraint} \label{sec:smoothsc}
Figure \ref{fig:bsqns} shows the PINNs solution to the Boussinesq equations \eqref{eq:Kookaburra} with fixed lambda $\lambda = 3$. The residue of each equation reaches to $10^{-4}$, except at the origin where the residue of each equation stays between 0.01 and 0.1, forming a peak. Systematic testing shows that even if more collocation points are chosen at the origin, the peak remains.  In fact, the peak is caused by the non-smoothness at the origin. We recall that the functional form of the neural network with activation function $\tanh$ is inherently a smooth function of the input. The non-smooth solution approximated by a PINN is expected to deviate from the true solution profile at the non-smooth point, resulting in large equation residues near the non-smooth point.  In other words, a unique way to reduce the error at the origin (non-smooth point) is to find the right $\lambda$ that corresponds to the smooth solution. Therefore, the constraint to find the smooth solution becomes equivalent to suppressing the error at the non-smooth point. Here, we impose an additional constraint using the higher-derivative of the residue of the equation around the origin as the {\it smoothness constraint}, namely
\begin{equation}\label{eq:loss_g}
	loss_{s}^{(k)}  = \frac{1}{N_s}\sum_{i=1}^{N_s} \left|\nabla f_{k} \left({\bf z}_i, \hat{q} \left({\bf z}_i \right) \right) \right|^2\,.
\end{equation}
where ${\bf z}_i$ indicates the random collocation points close to the origin (e.g. $|{\bf z}_i|<1$) and $N_s =500 $ is their total number.  $\hat{q}$ indicates the neural network prediction of each variable, and $f_k$ are the equation residues defined in \eqref{eq:bsqf}.  Then, the final cost function to discover the smooth solution to Bousinessq reads
\begin{equation}
J({\bf y, w}) = \tfrac{1}{n_b}\sum_{j=1}^{n_b} loss_{c}^{(j)} + \gamma \left( \tfrac{1}{n_e}\sum_{k=1}^{n_e} loss_{f}^{(k)} +\tfrac{1}{n_e}\sum_{k=1}^{n_e} loss_{s}^{(k)} \right)
\end{equation}
where $n_b=8$ and $n_e=6$ are the total number of solution conditions and governing equations as shown in \eqref{eq:bsqg} and \eqref{eq:bsqf}. For self-similar equations in this study, a good range of equation weight $\gamma$ to choose is between $O(0.01)$ and $O(0.1)$, within which the convergence rates of training are tested to be fast and relatively the same. However, we note that the training could fail when $\gamma > 1$. This is because all equations solved in this study are nonlinear equations, which involve many trivial solutions that are not of our interest. Proper boundary conditions need to be imposed to discover the right solution we want. However, when $\gamma>1$, the equation loss is considered a higher priority to be minimized during the training. Because of this, the training would be very easily trapped in the local minimum that corresponds to the trivial solution.

Figure \ref{fig:bsqsol} (same as Figure 1 in the main text) shows the smooth solutions to the Boussinesq equations derived by PINN with $\lambda = 1.917$ inferred by imposing the smoothness constraint \eqref{eq:loss_g}. The equation residues for the smooth solution are small and uniform everywhere in the domain of calculation. 

\subsection{Additional robustness tests of PINN solutions to Boussinesq equations}

\begin{figure}[t]
	\centering
	\includegraphics[width=0.8\textwidth]{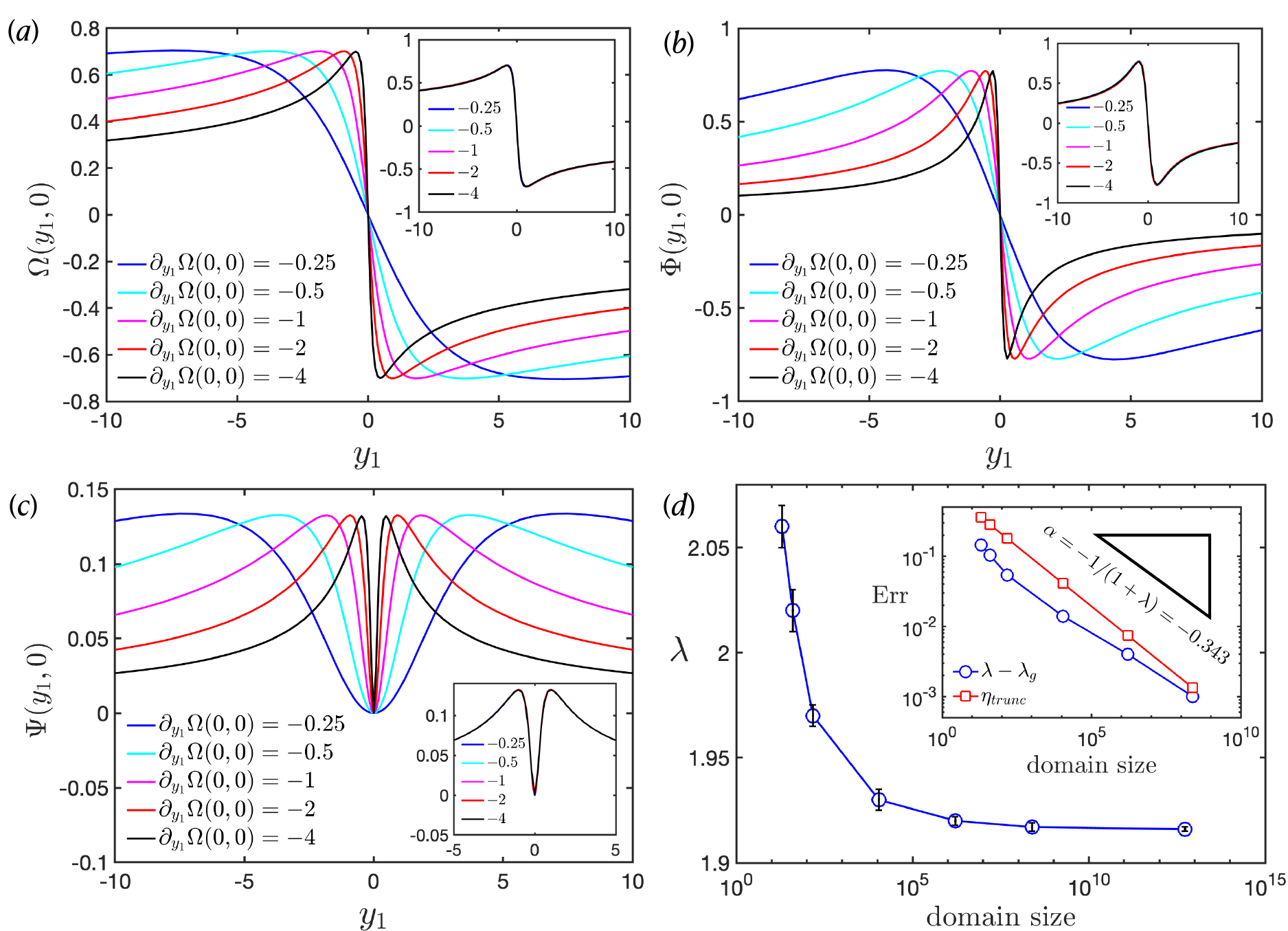}
	\caption{\label{fig:bsqanl} ($a$-$c$)  Solutions of $\Omega$, $\Phi$ and $\Psi$ (at $y_2=0$) derived from PINNs under different normalization conditions  show the same shape. The inset shows that all curves collapse on a single one after spatial rescaling.  ($d$) Inferred value for the smooth solution with respect to the domain size of calculation.  Error bars indicate the standard deviation of the inferred $\lambda$ from 4 independent PINN trainings with different random initialization. The inset shows the difference of the $\lambda$ inferred under smaller domain sizes with $\lambda_g$ which is inferred under the largest domain size here ($d =\textrm{sinh(30)} \approx 5\times10^{12}$). The difference between $\lambda$ is proportional to the truncation error $\eta_{trun}$ of the Hilbert transform due to the finite size of the domain.   }
\end{figure}

In order to ensure the robustness of the output of the PINN, we tested the performance of the method when the normalization condition is changed. In theory, the solution to the Boussinesq equation \eqref{eq:Kookaburra} has scaling symmetry. Namely, if $\Omega(y), \Psi(y), \Phi(y)$ is the solution to \eqref{eq:Kookaburra}, $\Omega\left(\frac{y}{\kappa}\right), \Psi\left(\frac{y}{\kappa}\right), \Phi\left(\frac{y}{\kappa}\right)$ for any $\kappa \neq 0$ are also solutions corresponding to a possibly different normalization.  Figure  \ref{fig:bsqanl} ($a$-$c$) shows the solution of $\Omega$, $\Phi$ and $\Psi$ in $y$-coordinates for different normalization conditions.  After rescaling by their corresponding $\kappa$, all solutions collapse well on a single curve (see inset),  confirming the scaling symmetry of the solution derived by PINNs.  The inferred $\lambda$ for different normalization conditions is consistent with $\lambda = 1.917 \pm 0.002$. Here, the error comes from experiments with different random initialization of independent training of PINN.

In addition, we also tested the stability of our method with respect to the domain for different domain sizes, observing the convergence of the value of $\lambda$ as the domain size grew towards infinity (Figure \ref{fig:bsqanl} $d$). Here,  we compute the error of $\lambda$ for different domain sizes by $\lambda-\lambda_g$, where $\lambda_g$ indicates the $\lambda$ inferred for the largest domain size in our test, $d_g = \sinh(30) \approx 5\times 10^{12}$.  The truncation error $\eta_{trunc}$ indicates the error of solutions due to the finite domain of calculation. Using formal asymptotics, one can obtain that $\Phi$ and $\Psi$ decay over $y_1$ by a power law with exponent $\alpha = -1/(1+\lambda)$,  which only reaches 0 when $y_1$ or $y_2$ goes to infinity. For a finite domain, imposing the constraint $\Phi(|y_1|=d)$ in the training, the truncation error of the approximate solution from PINN with the exact solution can be estimated by $\eta_{trunc} = d^{\alpha}$.  The inset of Figure \ref{fig:bsqanl}($c$) shows that the error of $\lambda$ in log-log scale is proportional to the truncation error $\eta_{trunc}$ of the finite domain size, which also decays over $d$ with a power law exponent $\alpha = -1/(1+\lambda)$.

\subsection{Comparison with literature}\label{ss:non-smooth}

Chen and Hou \cite{Chen2021} constructed self-similar solutions lying in the space $C^{1,\alpha}$ for $\alpha>0$ very small, where $\alpha = 1/(1+\lambda)$. The self-similar solutions are then used to prove blow-up for 3-D Euler in the presence of boundary. To describe Chen and Hou's result, we choose a large fixed $\lambda$, corresponding to a small $\alpha>0$,  and define
 \[R=\left(y_1^2+y_2^2\right)^{\frac{\alpha}{2}}\quad\mbox{and}\quad\gamma=\arctan\left(\frac{y_2}{y_1}\right)\,.\]
Chen and Hou \cite{Chen2021} showed that the self-similar Boussinesq equations \eqref{eq:Kookaburra} have an approximate analytic solution given by\footnote{We remark that Chen and Hou \cite{Chen2021} adopt the opposite sign convention for vorticity than that used in this paper.} \begin{gather}\label{eq:approx:sol}
 \Omega=-\frac{\alpha}{c}\left(\cos (\gamma)\right)^\alpha\frac{3R}{(1+R)^2},\quad
  \Phi=-\frac{\alpha}{c}\left(\cos (\gamma)\right)^\alpha\frac{6R}{(1+R)^3},\quad \Psi \equiv 0
 \end{gather}
 for $\gamma\in[0,\frac\pi 2]$ (or equivalently $y_1\geq0$) and where
 \[c=\frac{2}{\pi}\int_0^{\frac \pi 2}\left(\cos (\theta)\right)^\alpha\sin(2\theta)\,d\theta\,.\]
 One then extends the approximate solution $ \Omega$ and  $ \Phi$ to the region $y_1\leq 0$ via an odd extension.
\begin{figure}[t]
	\centering
	\includegraphics[width=\textwidth]{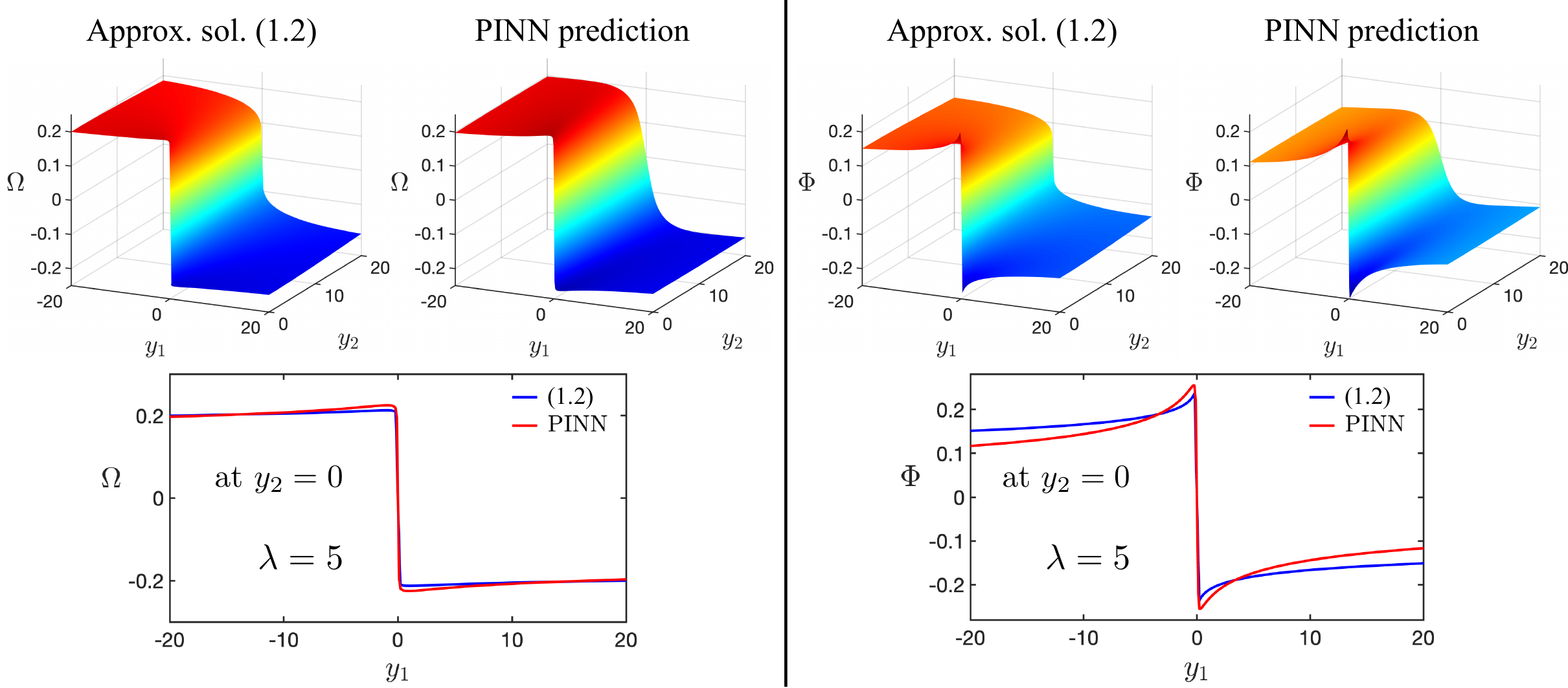}
	\caption{\label{fig:nonsm} Non-smooth solutions of $\Omega$ and $\Phi$ for the 2D Bousinessq equations \eqref{eq:Kookaburra} derived by the physics-informed neural network for $\lambda = 5$, which shows an agreement with the approximate solution \eqref{eq:approx:sol}. }
\end{figure}

Figure \ref{fig:nonsm} shows that for $\lambda=5$, the solution found by PINNs agrees well with the approximate solution \eqref{eq:approx:sol}.\footnote{We stress that the solution \eqref{eq:approx:sol} is a first order approximation, whereas the PINN is attempting to approximate the exact solution. In particular, a discrepancy between the two solutions is expected.}  In addition to providing an opportunity to validate the numerical method, the ability to find non-smooth solutions opens the possibility of using a continuation argument to search for smooth solutions. In the end, the PINN method is sufficiently robust that no such continuation argument is necessary; however, the possibility of such a continuation argument may be advantageous for other settings, such as finding unstable solutions as mentioned at the end of the main manuscript.

\section{The Burgers' equation}
\subsection{Problem setting}
One of the simplest PDEs which exhibits self-similar blow-up is the 1-D Burgers' equation:
\begin{equation*}
u_t+uu_x=0\,.
\end{equation*}
Assuming the self-similar ansatz
 \begin{equation*}
u=(1-t)^{\lambda}U\left(\tfrac{x}{(1-t)^{1+\lambda}}\right)\,,
 \end{equation*}
and defining $y=\tfrac{x}{(1-t)^{1+\lambda}}$,
 we obtain the self-similar Burgers' equation
  \begin{equation}\label{eq:SS:Burgers}
-\lambda U+((1+\lambda)y+ U)\partial_y U=0\,.
\end{equation}
We impose $U$ to be odd. The self-similar Burgers' equation can then be solved implicitly by
\begin{equation}\label{eq:implicit:sol}
y=-U-CU^{1+\frac 1\lambda}
\end{equation}
for some constant $C$. To determine $C$, we add the constraint
\begin{equation*}
U(-2)=1\,,
\end{equation*}
which fixes $C=1$, independent of $\lambda$. It is classical \cite{EgFo2009}, that in order to obtain a smooth self-similar solution, then $\lambda$ must be chosen such that\begin{equation*}
\lambda = \frac{1}{2i+2}\quad\mbox{for }i=0,1,2,\dots
\end{equation*}
Otherwise, the solution \eqref{eq:implicit:sol} has non-smoothness at $y=0$ (origin).  Due to the availability of the analytic solution,  Burgers' equation becomes a good touchstone to test our method. 
\begin{figure}[t]
	\centering
	\includegraphics[width=0.85\textwidth]{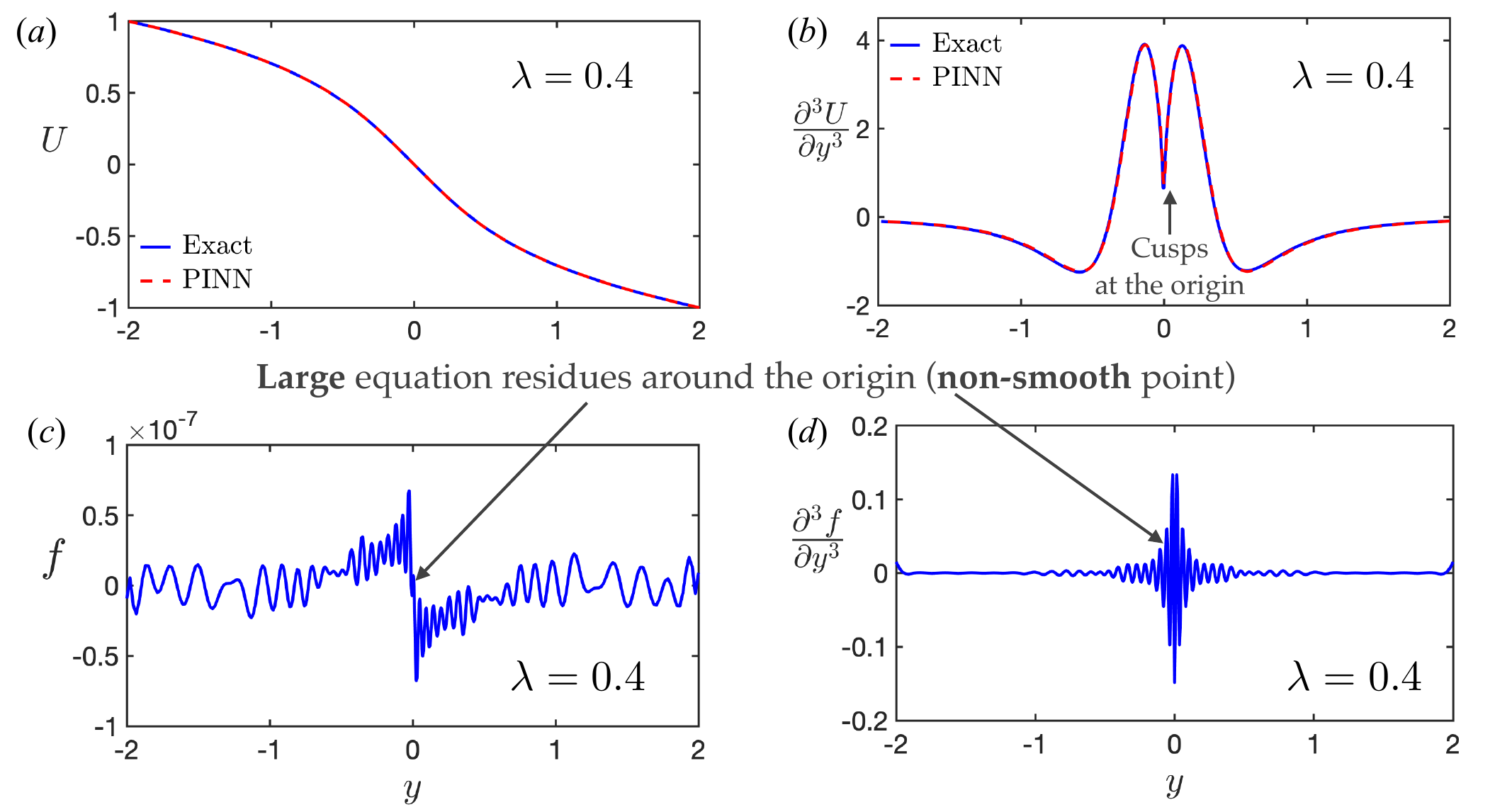} \vspace{-1em}
	\caption{\label{fig:bgns} ($a$) A non-smooth solution for the Burgers' equation with fixed $\lambda = 0.4$ using a physics-informed neural network. ($b$) PINN successfully captures the cusp at the third derivatives of the solution. ($c$) Equation residue and ($d$) its high-order derivative at the non-smooth point ($y=0$) are larger than the rest of domain.  }
\end{figure}

\begin{figure}[t]
	\centering
	\includegraphics[width=\textwidth]{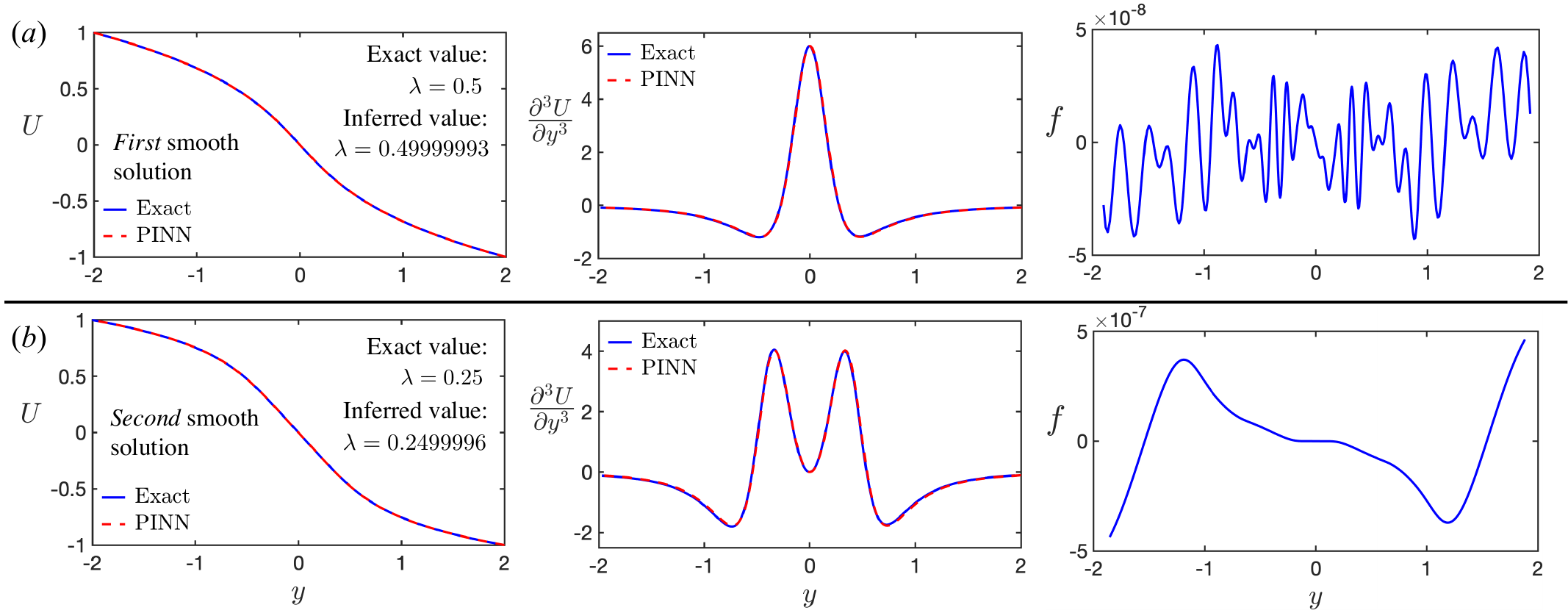}\vspace{-1em}
	\caption{\label{fig:burgers} The first and second smooth solutions for the Burgers' equation derived from the physics-informed neural network. The error of the inferred $\lambda$ for the first smooth solution is of order $10^{-8}$ and the second of order $10^{-7}$, which is consistent with the equation residue. }
\end{figure}

\subsection{PINN setting} 
We impose odd symmetry of the variable $U$ by constructing the function form $U = [\mathrm{NN}_u(y) - \mathrm{NN}_u(-y)]/2$, where $\mathrm{NN}_u$ indicates a multi-layer fully-connected neural network created for $U$. Here, we use 3 hidden layers with 20 units in each hidden layer for the network  $\mathrm{NN}_u$ and use the hyperbolic tangent function $\tanh$ as the activation function. The {\it condition loss} and {\it equation loss} for the Burgers' equation are
\begin{equation}\label{eq:bgf}
	g  = U(-2) - 1    \qquad \text{and} \qquad  f(y, U) = -\lambda U+((1+\lambda)y+ U)\partial_y U
\end{equation} 

\subsection{Non-smooth solution}
By fixing $\lambda = 0.4$, based on the analytic expression \eqref{eq:implicit:sol},  the solution has no smoothness at the origin. Figure \ref{fig:bgns} shows the solution of $U$ derived by PINN, which is in good agreement with the analytic solution. In addition, the solution from PINN can capture the cusp at the origin (non-smoothness point) in the third derivative of the solution, indicating the accuracy of the method.  However, as mentioned in \S \ref{sec:smoothsc}, the neural network with $\tanh$ activation function is inherently a smooth function, which is expected to deviate more from the non-smooth solution at the non-smooth point. This is confirmed by the equation residue $f$ (Figure \ref{fig:bgns}$c$), which shows a larger error at the origin (non-smooth point) than anywhere else.  The error at the non-smooth point is further amplified at the higher-order derivative of the residue $f$ (Figure \ref{fig:bgns}$d$).



\subsection{Smoothness constraint and smooth solutions}
The explicit solution \eqref{eq:implicit:sol} to Burgers' equation shows that there exist more than one values of $\lambda$ leading to smooth solutions. In order to find a specific smooth self-similar solution, we constrain the parameter $\lambda$ to a window $\lambda\in[\frac{1}{2i+3},\frac{1}{2i+1}]$.  In addition, to single out the smooth solution at $\lambda = 0.5$ from the non-smooth solution nearby that $\lambda$, we penalize the third derivative of the equation residue. Here, we are utilizing the fact that the non-smooth solutions in the neighborhood of $\lambda = 0.5$ has unbounded fourth-order derivative, which appears in the third order of derivative of equation residue.
 Namely, the smoothness constraint for the first smooth solution is
\begin{equation}\label{eq:loss_gb1}
	loss_{s1}  = \frac{1}{N_s}\sum_{i=1}^{N_s} \left|\frac{d^3 f}{d y^3} \left(y_i, \: U(y_i ) \right) \right|^2\,.
\end{equation}
where $y_i$ indicates the random collocation points close to the origin (e.g. $|y_i|<1$) and $N_s$ is their total number. $f(y, U)$ is the equation residue defined in \eqref{eq:bgf}.  However, for the second smooth solution $\lambda=0.25$, the non-smooth solutions near $\lambda = 0.25$ have unbounded sixth-order derivative. In order to distinguish the second smooth solution from the non-smooth solution, we need to impose a smoothness constraint by adding the {\it fifth} derivative of the equation residue around the origin, namely,
\begin{equation}\label{eq:loss_gb1}
	loss_{s2}  = \frac{1}{N_s}\sum_{i=1}^{N_s} \left|\frac{d^5 f}{d y^5} \left(y_i, \: U(y_i ) \right) \right|^2\,.
\end{equation}
Consequently,  further higher-order derivatives of equation residue $f$ are required to find more smooth solutions at smaller $\lambda$, which would take a larger amount of computational time for the training. 
By imposing the right smoothness constraint, Figure \ref{fig:burgers} demonstrates that the PINN can find the smooth solutions for the self-similar Burgers equation \eqref{eq:SS:Burgers} along with their corresponding self-similar parameter $\lambda$ with high precision.  In fact, the case $\lambda=0.5$ ($i=0$) corresponds to stable, smooth self-similar profiles, while the cases $\lambda=0.25$ ($i=1$) correspond to an unstable self-similar profile. \color{black}

\section{The generalized De Gregorio equation}

\subsection{Problem setting}
Consider the generalized De Gregorio equation  \cite{OkSaWu08}
 \begin{equation*}
 \omega_t+au\omega_x=\omega u_x,\qquad\mbox{where } u=\int_0^x(H\omega)(s)\,ds
 \end{equation*}
and $H$ is the Hilbert transform.  We assume the following self-similar ansatz
 \begin{equation*}
 \omega=\frac{1}{1-t}\Omega\left(\frac{x}{(1-t)^{1+\lambda}}\right)\,.
 \end{equation*}
 Then, if we define $U=\Lambda^{-1}\Omega$, make the change of coordinates  $y=\tfrac{x}{(1-t)^{1+\lambda}}$,  we obtain the self-similar equations
 \begin{equation}\label{eq:deg}
 	 \Omega+((1+\lambda)y+a U)\partial_y \Omega-\Omega \partial_y U=0 \quad \text{where }  U = \int_{0}^{y} H\Omega(s)\:\mathrm{d}s
 \end{equation}
Note that the equation is invariant under the rescaling $\Omega (y)\mapsto \Omega(\mu y)$. As such, we are free to fix
 \begin{equation}\label{e:cond-1}
 \partial_y \Omega(0)= -2\ .
 \end{equation}
 We assume that both $\Omega$ and $U$ are odd. The generalized De Gregorio equation in the case of $a=-1$ is also known as the C\'ordoba-C\'ordoba-Fontelos (CCF) equation.

\subsection{PINN setting}
We impose odd symmetry of $U$ and $\Omega$ by constructing the function form $q = [\mathrm{NN}_q(y) - \mathrm{NN}_q(-y)]/2$, where $\mathrm{NN}_q$ indicates a multi-layer fully-connected neural network created for each variable $q$. Here, We use 3 hidden layers with 20 units for each network  $\mathrm{NN}_q$ and use the hyperbolic tangent function $\tanh$ as the activation function. The {\it condition loss} for the De Gregorio equation is
\begin{equation}
	g_1  = \partial_y\Omega(0) + 2 , \qquad 	g_2 = \nabla  {\bf U}(|y|=d),  \quad \text{and} \quad 	g_3=\Omega(|y|=d) \\
\end{equation} 
where $d$ is the domain size of the calculation.  Due to the nonlocality of the Hilbert Transform, the high-accuracy solution to the generalized De Gregorio equation should be solved in a very large domain. To cope with that, we defined a new $z$-coordinate, similar to that for the Boussinesq equations,  with relation
\[ y= \sinh(z)  \qquad \Longleftrightarrow  \qquad z =\sinh^{-1} y \,.\] 
The self-similar generalized De Gregorio equations in $z$-coordinates can then be written as:
\begin{equation}
\Omega+((1+\lambda)\sinh(z)+a U)\sech(z)\partial_z \Omega-\Omega \sech(z)\partial_z U=0 \qquad \text{with } \quad	\sech{z}\frac{\partial U}{\partial z}= H\Omega(z).
\end{equation}
The latter equation is derived from $\partial_y U = H\Omega(y)$, which is equivalent to that in \eqref{eq:deg}, considering that the odd symmetry of $U$ yields the implicit boundary condition $U(0) = 0$. Thus,  the {\it equation loss} contains the residue of two equations:
\begin{equation}\label{eq:degs}
	f_1  = \Omega+((1+\lambda)\sinh(z)+a U)\sech(z)\partial_z \Omega-\Omega \sech(z)\partial_z U  \quad \text{and} \quad  f_2 = \sech{z}\frac{\partial U}{\partial z} -\tilde{H}\Omega(z) \\
\end{equation} 
where $\tilde{H}$ indicates the numerical calculation of the Hilbert Transform that is described in the next section.

  \begin{figure}[t]
	\centering
	\includegraphics[width=0.9\textwidth]{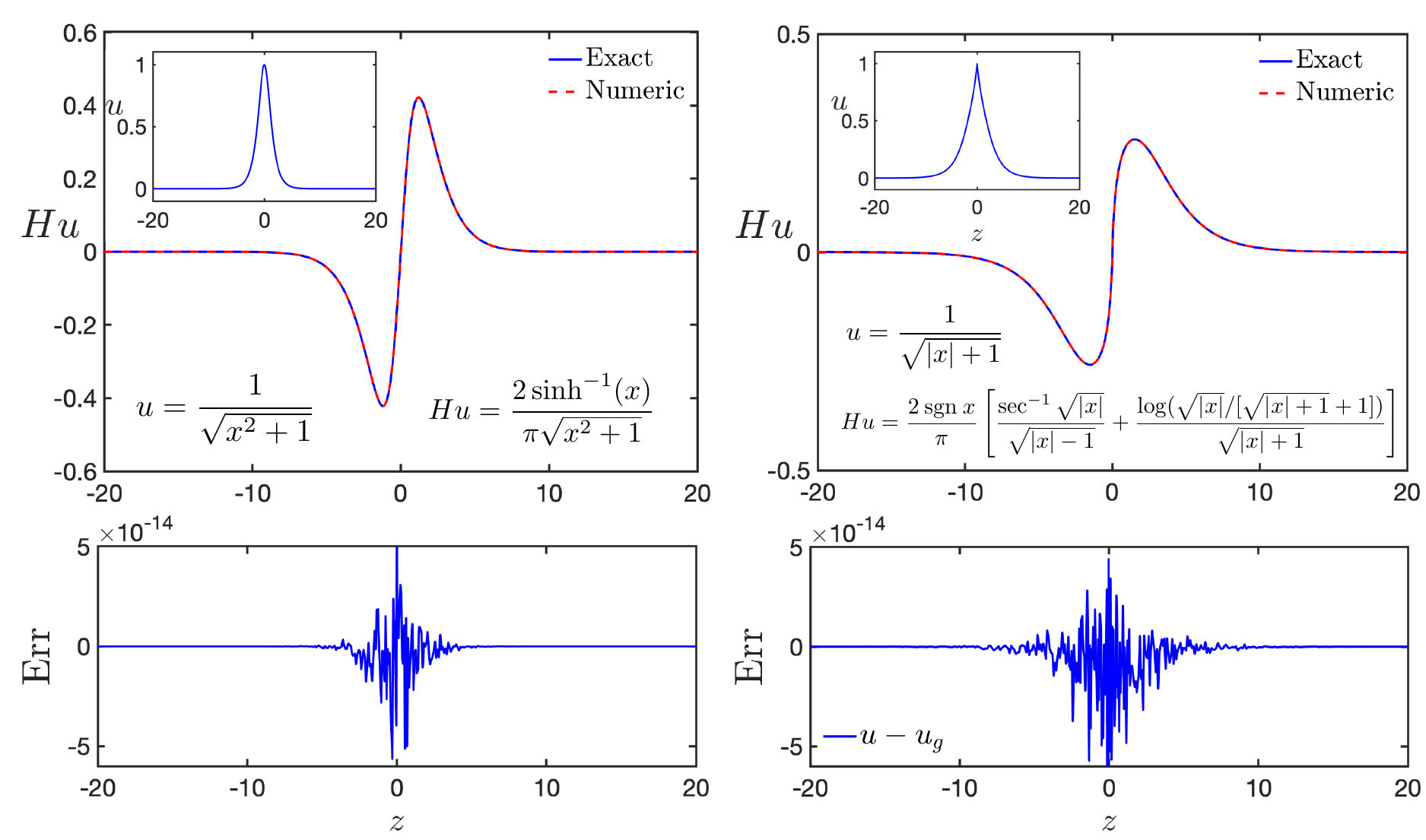} \vspace{-1em}
	\caption{\label{fig:ht_acc} ($a$, $c$)  Comparison of numerical Hilbert transform developed here with the exact expressions for two slow-decaying or non-smooth functions in the $\sinh^{-1}$ coordinates, where $z = \sinh^{-1}(y)$. By using 50000 sampling points with equal distance in the $z$-coordinate, the error of the numerical Hilbert transform reaches the round-off error of double precision, confirming the accuracy of the method.}
\end{figure}

\subsection{Implementation of the Hilbert Transform}  \vspace{-0.2em}
The difficulty of solving the generalized De Gregorio equation lies in the presence of the Hilbert Transform.  In order to perform an accurate calculation of the nonlocality given by the Hilbert Transform, and in particular due to the very large domains,  we extended the method of Bilato--Maj--Brambilla \cite{bilato2014algorithm} by using third-order Lagrange interpolations to compute the Hilbert Transform integral in each local interval and generalize the method to the $z$-coordinates, which significantly reduce the truncation error. Figure \ref{fig:ht_acc} shows the comparison of the numerical Hilbert transform using our method with the analytic Hilbert transform for two slow-decaying functions. By using 50,000 sampling points equally-spaced in the $z$-coordinate, our new method yields very small errors (both truncation and integration errors are of the order of machine-size error).  Here, the sampling point refers to the position where we evaluate the integrand of the Hilbert Transform. In addition, our method is able to calculate the Hilbert Transform of non-smooth functions (Figure \ref{fig:ht_acc}$b$) with the same high precision. Compared with that, the classical FFT-based methods yield insufficient precision due to the slow-decaying tails of the functions. Here, we consistently use our newly-developed numerical Hilbert Transform method to calculate the residue \eqref{eq:degs} of the generalized De Gregorio equation.

 \begin{figure}[t]
	\centering
	\includegraphics[width=\textwidth]{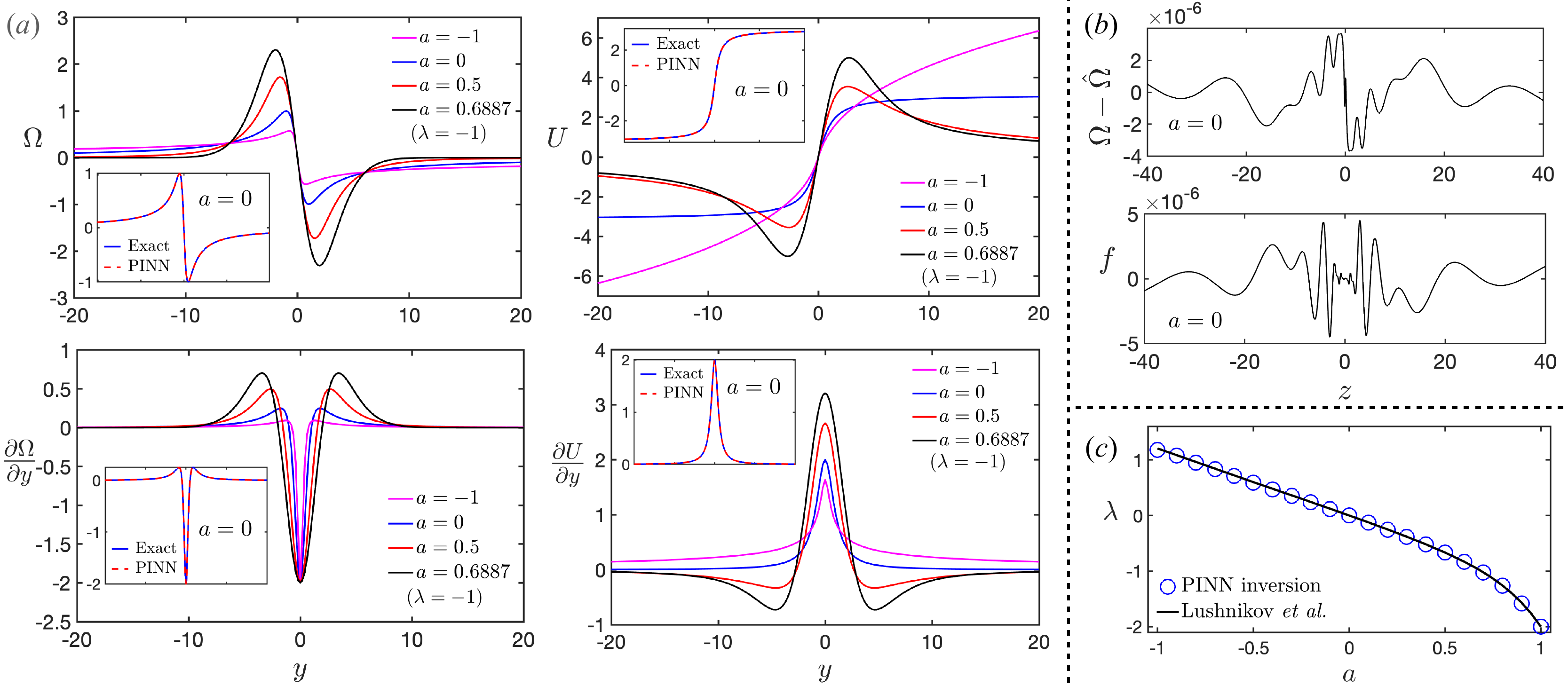} \vspace{-2em}
	\caption{\label{fig:DeGreg} ($a$) Smooth solution $\Omega$ and $U$ to the generalized De Gregorio equation for different values of $a$. Besides inferring $\lambda$ for a given value of $a$, PINN can also directly infer $a$ for a given value of $\lambda$. The black line indicates the solution inferred by imposing $\lambda = -1$, which gives $a = 0.6887$, which shows a good agreement with the result from \cite{Lushnikov2021}.  The inset shows the comparison  of the PINN prediction with the exact solution for $a=0$, which shows a perfect agreement.  ($b$) The error between the PINN solution to the  Constantin-Lax-Majda equation (generalized De Gregorio equation with $a=0$) and the exact solution \eqref{eq:CLM}, which is of the same order of magnitude as that of the equation residue.  ($c$) the curve of $\lambda$ inferred via PINNs for different values of $a$ within the range of $[-1, 1]$ .  The curve collapses well with the result from \cite{Lushnikov2021}, indicating the accuracy of PINN inversion. }
\end{figure}

\subsection{Smooth solution and its comparison with literature} \vspace{-0.2em}
By imposing the smoothness constraint that penalizes the derivative of equation residue \eqref{eq:degs} around the origin (non-smooth point), we can infer the $\lambda$ that corresponds to the smooth solution for each value of $a$ between $[-1, 1]$ (Figure \ref{fig:DeGreg}),  which reproduce the findings of \cite{Lushnikov2021} including the prediction that self-similar collapse forms for $a$ below the critical value $a_c\approx0.6887$.  When $a=0$, the generalized De Gregorio equation turns into the Constantin-Lax-Majda equation, which has the following analytic solution with respect to the normalization condition \eqref{e:cond-1}:
\begin{equation} \label{eq:CLM}
	\hat{\Omega}(y) = -\frac{2y}{1+y^2} \qquad \text{and} \qquad \hat{U}(y) = 2\arctan (y).
\end{equation}
The inset of Figure \ref{fig:DeGreg} shows that the solution from PINN has a good agreement with the analytic solution \eqref{eq:CLM}. The relative error between them is of order $10^{-6}$  (Figure \ref{fig:DeGreg}$b$), consistent with that of the equation residue.

\subsection{Smooth solution and smooth constraint to the CCF equations}
When $a=-1$, the general De Gregorio equation is known as the C\'ordoba-C\'ordoba-Fontelos (CCF) equation. The first smooth solution as shown in Figure \ref{fig:DeGreg} is found by adding the smoothness constraint 
\begin{equation}\label{eq:loss_ccfs1}
	loss_{s1}^{(k)}  = \frac{1}{N_s}\sum_{i=1}^{N_s} \left|\frac{d f_k}{d y} \left(z_i, \: \hat{q}(z_i ) \right) \right|^2\, ,
\end{equation}
where $z_i$ indicates the random collocation points close to the origin (e.g.\ $|z_i|<1$) and $N_s=80$ is their total number.  $\hat{q}$ indicates the neural network prediction of $U$ and $\Omega$, and $f_k(z, \hat{q})$ indicates the equation residues defined in \eqref{eq:degs}. As described in the main text, an open question for the CCF equation is whether there exist blow-up solutions with $\lambda < 1$. 

 \begin{figure}[t]
	\centering
	\includegraphics[width=0.70	\textwidth]{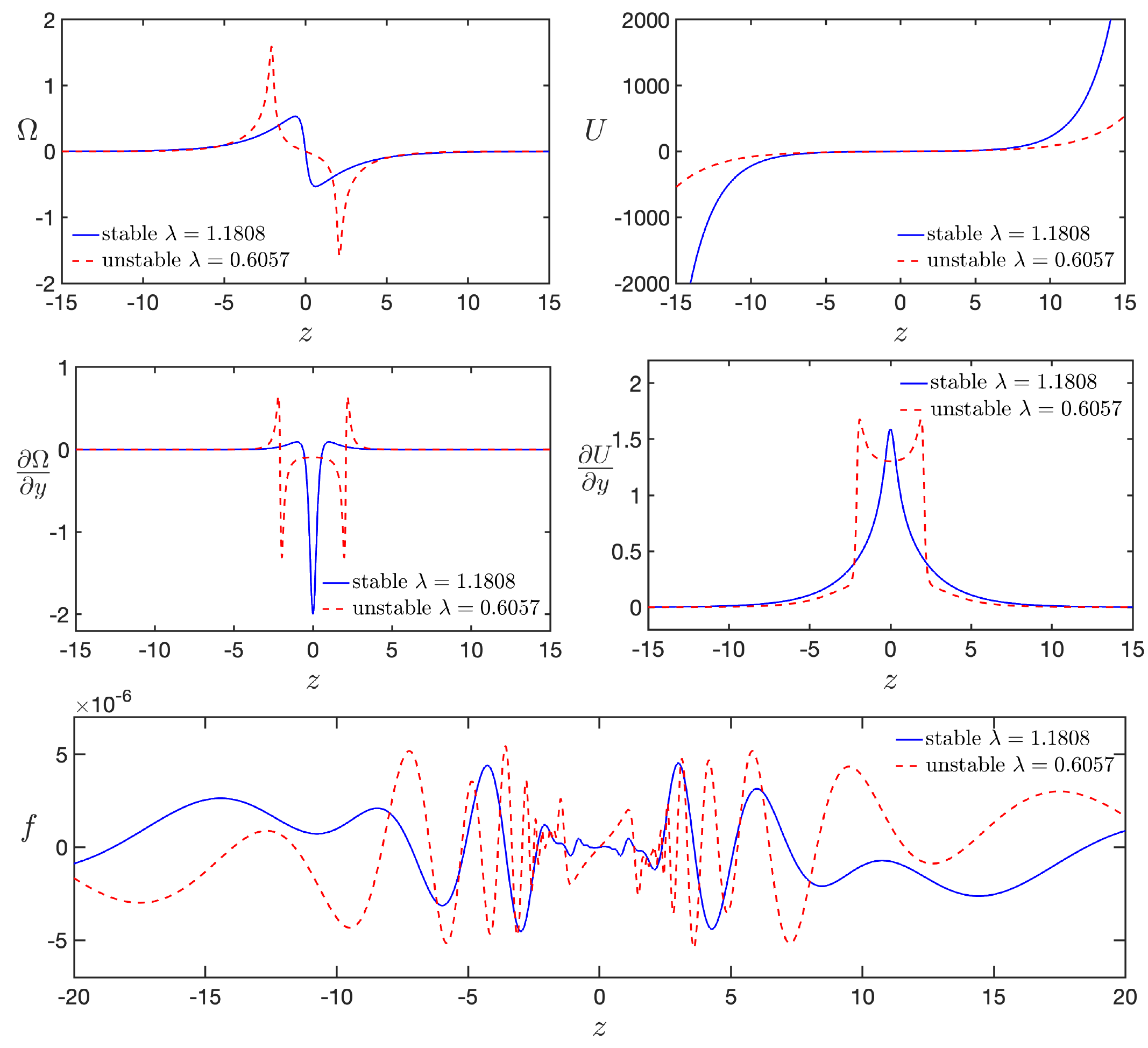} \vspace{-1em}
	\caption{\label{fig:ccf_s} ($a$-$d$) First (stable) and second (unstable) smooth solutions of the CCF equation in the $\sinh^{-1}$ coordinates. ($e$) The equation residues for both solutions reach $10^{-6}$, indicating their accuracy. }
\end{figure}

As shown in Figure \ref{fig:DeGreg}, the first smooth solution to the CCF equation is found at $\lambda = 1.1808$, which does not satisfy this criterion. However, analogous to the Burgers equation, there should also exist multiple smooth solutions to the self-similar CCF equation with different values of $\lambda$.  Learning from the second smooth solution to the Burgers' equation, we expect that the non-smooth solutions near the second smooth solution to the CCF at the smaller $\lambda$ should also have unbounded higher-order derivatives. Thus,  to distinguish the second smooth solution from the non-smooth ones, we impose a smoothness constraint by adding the {\it third} derivative of the equation residue around the origin, namely,
\begin{equation}\label{eq:loss_ccfs2}
	loss_{s2}^{(k)}  = \frac{1}{N_s}\sum_{i=1}^{N_s} \left|\frac{d^3 f_k}{d y^3} \left(y_i, \: \hat{q}(y_i ) \right) \right|^2\,.
\end{equation}
By imposing the new smoothness constraint \eqref{eq:loss_ccfs2} and the new normalization condition $\Omega(0.5) = -0.05$, we successfully discovered, for the first time, the second smooth solution to the CCF equation with \mbox{$\lambda = 0.60573$}, which is also the first smooth solution discovered with $\lambda < 1$. Figure \ref{fig:ccf_s} shows the comparison of the shape between the first and second smooth solutions to the CCF. The equation residue for both smooth solutions can reach $10^{-6}$, confirming the accuracy of our PINN training. More robustness tests and stability analysis for both smooth solutions are  given in the next section.

\subsection{Additional robustness tests of the PINN solutions to the CCF equations}
\begin{figure}[t]
	\centering
	\includegraphics[width=0.85\textwidth]{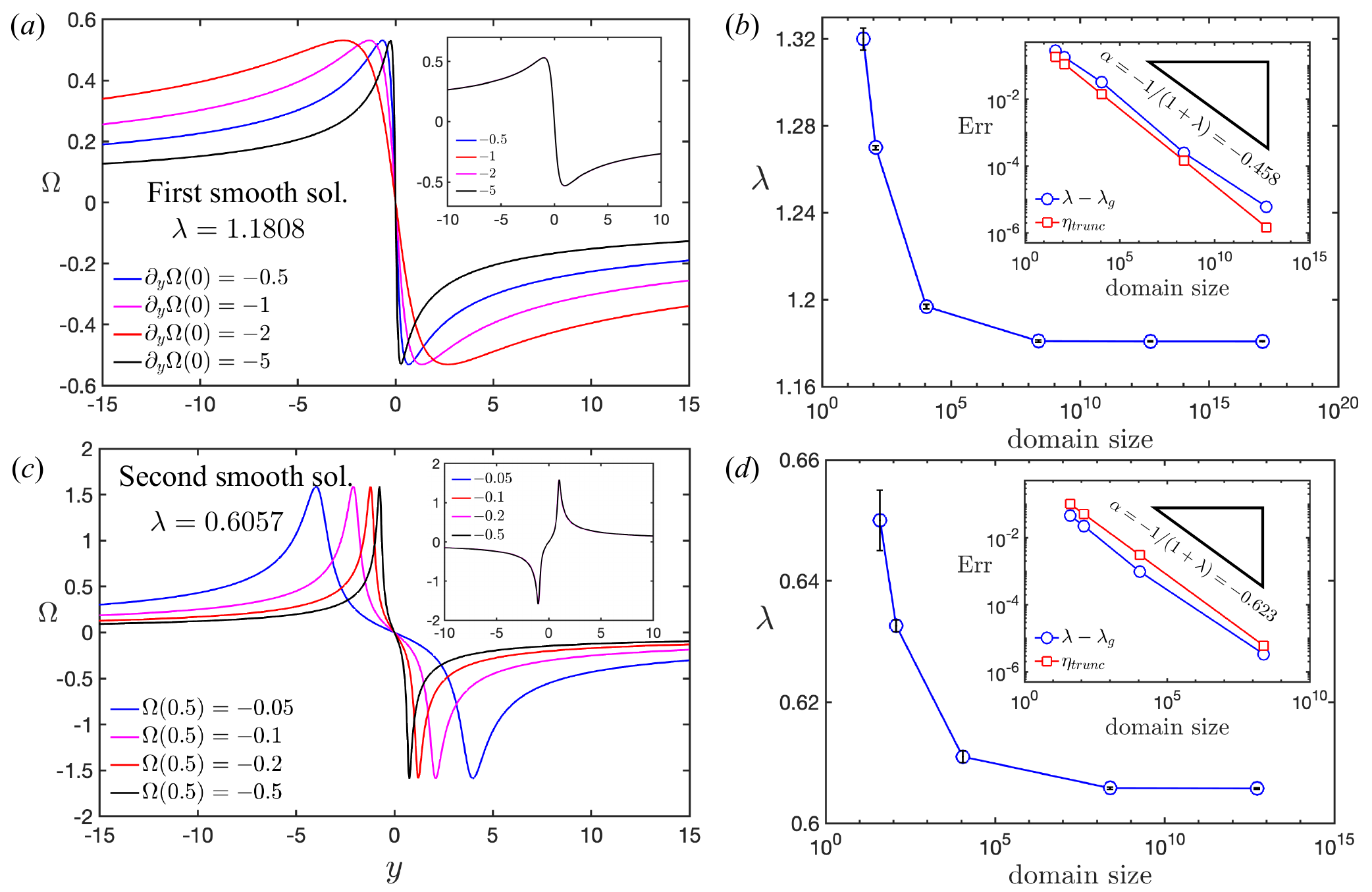} \vspace{-1em}
	\caption{\label{fig:ccf_anl} ($a$, $c$)  Solution of $\Omega$ from PINNs under different normalization conditions shows the same shape. The inset shows that all curves collapse on a single one after spatial rescaling.  ($b$, $d$) Inferred value for the smooth solution with respect to the domain size of calculation.  Error bars indicate the standard deviation of the inferred $\lambda$ from 3 independent PINN trainings with different random initialization. The inset shows the difference of the $\lambda$ inferred under smaller domain sizes with $\lambda_g$ inferred under the largest domain size here ($d =\textrm{sinh(40)} \approx 10^{17}$ for stable and $d =\textrm{sinh(30)} \approx 5\times10^{12}$ for unstable). The difference between $\lambda$ is proportional to the truncation error $\eta_{trun}$ of the Hilbert Transform due to the finite size of the domain.}
\end{figure}

As in the robustness test of PINN solution to the Boussinesq equation,  we tested the performance of PINN method on the solution to the CCF equation when the normalization condition is changed. The solution to the CCF equation \eqref{eq:deg} has scaling symmetry. Namely, if $\Omega(y)$ is the solution to \eqref{eq:deg}, $\Omega\left(\frac{y}{\kappa}\right)$ for any $\kappa \neq 0$ are also solutions, but may correspond to different normalization conditions.  Figure  \ref{fig:ccf_anl} ($a$, $c$) shows the first and second smooth solutions of $\Omega$ in the $y$-coordinates for different normalization conditions.  After rescaling with their corresponding $\kappa$, all solutions collapse well on a single curve (see inset),  confirming the scaling symmetry of the solution derived by PINNs.  The inferred $\lambda$ for different normalization conditions is consistent with $\lambda = 1.1808 \pm 0.00005$ and  $\lambda = 0.60573 \pm 0.00004$ for the first and second smooth solutions, respectively. 

In addition, we also tested the stability of our method with respect to the domain for different domain sizes, observing the convergence of the value of $\lambda$ for both the first and second smooth solutions  (Figure \ref{fig:ccf_anl} $b$, $d$) as the domain size grew towards infinity, the same as shown in Figure \ref{fig:bsqanl}($d$).  Here,  we compute the error of $\lambda$ for different domain sizes by $\lambda-\lambda_g$, where $\lambda_g$ indicates the $\lambda$ inferred for the largest domain size in our test, e.g. $d_g = \sinh(40) = 1.2\times 10^{17}$ for the first smooth solution and $d_g =\sinh(30) = 5\times 10^{12}$ for the second smooth solution.  The truncation error $\eta_{trunc}$ indicates the error of the Hilbert Transform due to the finite domain of calculation. The solution of $\Omega$ decays over $y$ by a power law with the exponent $\alpha = -1/(1+\lambda)$,  which reaches 0 only when $y$ goes to infinity. However, by calculating the Hilbert transform of $\Omega$ within the finite domain within $[-d, d]$ in the training, the truncation error of the approximate solution from PINN with the exact solution can be estimated by $\eta_{trunc} \sim d^{\alpha}$.  The inset of Figure \ref{fig:ccf_anl}($b$, $d$) shows that the error of $\lambda$ in log-log scale is proportional to the truncation error $\eta_{trunc}$ of the finite domain size, which also decays over $d$ with a power law exponent $d = -1/(1+\lambda)$.

\begin{figure}[t]
	\centering
	\includegraphics[width=1\textwidth]{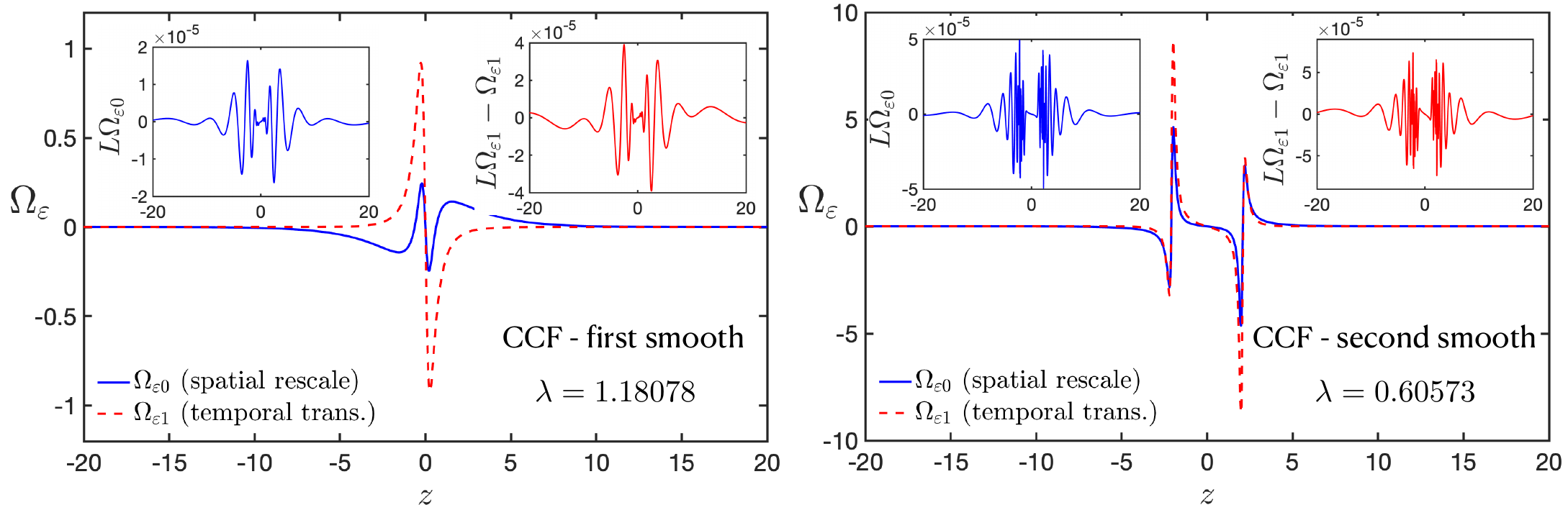} \vspace{-1.5em}
	\caption{\label{fig:ccf_eig}  Eigenfunctions $\Omega_{\epsilon 0}$ and $\Omega_{\epsilon 1}$ for the first and second smooth solutions to the CCF equation in the $z$-coordinate derived from PINNs, which correspond to the spatial rescaling and temporal translation properties of the solution. The inset shows that the analytic expression of both trivial eigenfunctions for both smooth solutions satisfy the linearized equation, which confirms the accuracy of the solution.}
\end{figure}

\subsection{Stability of the second smooth solution to the CCF equation}
\begin{figure}[t]
	\centering
	\includegraphics[width=1\textwidth]{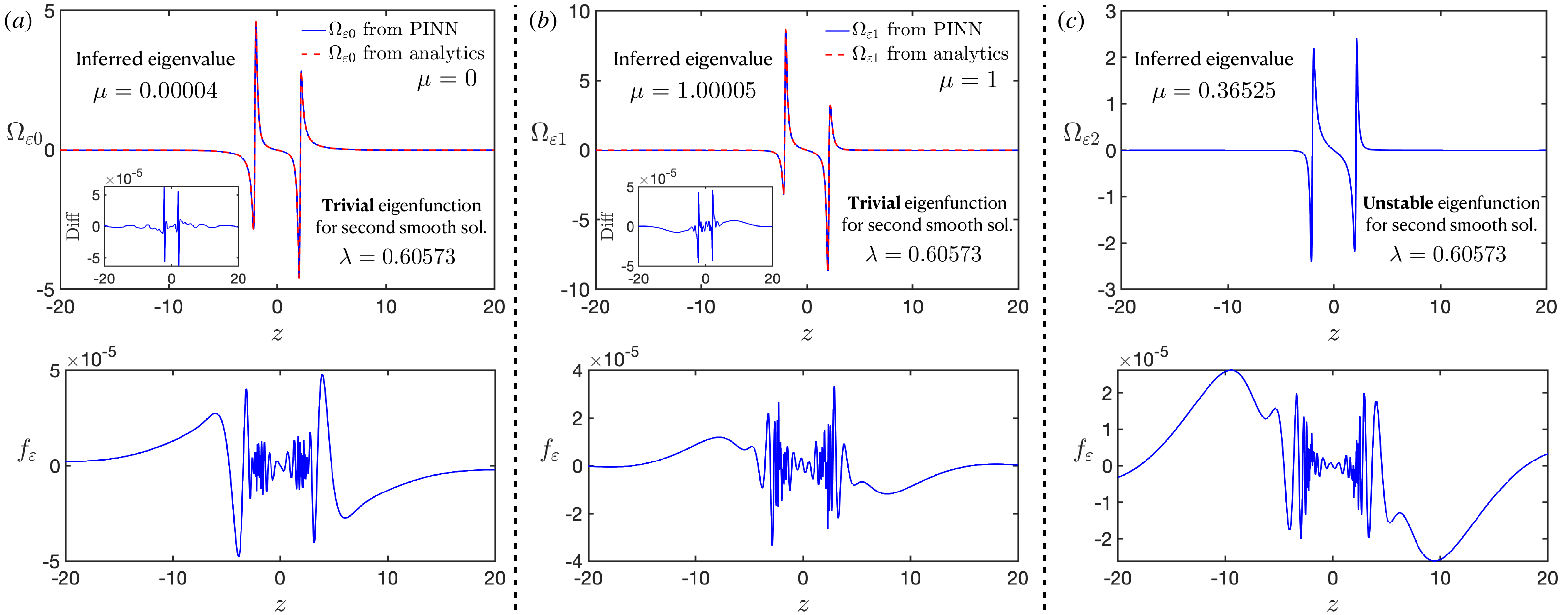} \vspace{-1.5em}
	\caption{\label{fig:ccf_eig2}  Eigenfunction $\Omega_\epsilon$ with associated eigenvalue $\mu$ for the second smooth solution to the CCF equation derived from PINN by solving the linearized equation \eqref{eq:lineareqn} considering free parameter $\mu$.  ($a$, $b$) For the two trivial eigenvalue/eigenfunction  the PINN results shows a good agreement with their analytic expressions correspond to the spatial rescaling \eqref{eq:eigspt} and temporal translation \eqref{eq:eigtp} properties of the solution. The difference between the PINN and analytics is consistent with the equation residue $f_\epsilon$ of the linarized equation, namely $f_\epsilon = L\Omega_\epsilon - \mu \Omega_\epsilon$.  ($c$) In addition, we found the third eigenfunction associated with eigenvalue $\mu = 0.36525>0$, which confirm the second smooth solution to the CCF equation is an {\it unstable} solution.  }
\end{figure}

The second smooth solution for the Burgers' equation at $\lambda = 0.25$ is an unstable solution. To determine whether the second smooth solution to the CCF equation at $\lambda = 0.6057$ is stable or unstable, we conduct the stability analysis. From now on, we suppose $(\hat\lambda,\hat \Omega,\hat U)$ is an approximate solution to the self similar CCF equation derived from PINN. By introducing the perturbation terms
\begin{equation}
	\Omega = \hat\Omega + \Omega_\epsilon \qquad  \text{and}  \qquad U = \hat U + U_\epsilon,
\end{equation}
we can define the linearization of the CCF equation around our approximate solution as
\begin{equation}\label{eq:linear:op}
	\mathcal L \Omega_\epsilon= -\Omega_\epsilon -((1+\hat\lambda)y- \hat U)\partial_y \Omega_\epsilon
	+U_\epsilon\partial_y \hat \Omega+\hat\Omega H\Omega_\epsilon +\Omega_\epsilon \partial_y \hat U  \quad \text{with} \ U_\epsilon=\int_0^y H\Omega_\epsilon\,.
\end{equation}
To calculate the eigenvalue and associated eigenfunction associated with each smooth solution to the self-similar CCF equation, we look for solutions to the equation
\begin{equation}\label{eq:lineareqn}
	\mathcal L(\Omega_\epsilon, U_\epsilon)=\mu \Omega_\epsilon
\end{equation}
where $\mu\in\mathbb C$ is the eigenvalue of the solution that needs to be determined. For each smooth solution to CCF, there exist two trivial eigenvalues/eigenfunctions corresponding to scaling, which are 
\begin{enumerate}
	\item  {\bf spatial rescaling}: $\mu = 0$ with corresponding eigenfunction 
	\begin{equation}\label{eq:eigspt}
		\Omega_{\epsilon 0} =	y\partial _y \hat\Omega
	\end{equation}
	\item {\bf temporal translation}: $\mu = 1$ with corresponding eigenfunction 
	\begin{equation}\label{eq:eigtp}
		\Omega_{\epsilon 1} =\hat\Omega+(1+\hat \lambda)y \partial_y \hat \Omega
	\end{equation}
\end{enumerate}  
Using the approximate solution from PINN with the theoretical expression of the two trivial eigenfunctions given above,  Figure \ref{fig:ccf_eig} shows that the two eigenfunctions $\Omega_{\epsilon 0}, \Omega_{\epsilon 1}$ of both the first and second smooth solutions, do satisfy the linearized equations with high precision,
\begin{equation}
	\mathcal{L} \Omega_{\epsilon 0} = 0\cdot \Omega_{\epsilon 0}, \qquad \text{and} \qquad \mathcal{L} \Omega_{\epsilon 1} = 1 \cdot \Omega_{\epsilon 1} 
\end{equation}
This assessment further confirms the accuracy of the two smooth solutions derived from PINNs, which contains the scaling property of the exact solution to the equation.  Besides the above assessment, we can also solve the linearized equation \eqref{eq:lineareqn} of $\Omega_{\epsilon}$ and $U_\epsilon$ considering free parameter $\mu$ via PINNs.  Figure \ref{fig:ccf_eig2}($a$, $b$) shows the two trivial eigenvalues and eigenfunctions discovered by PINN for the second smooth solution at $\lambda = 0.6057$, which shows a good agreement with the analytic eigenvalue and eigenfunction corresponding to the spatial rescaling \eqref{eq:eigspt} and temporal translation \eqref{eq:eigtp} properties of the solution.

To confirm whether the second smooth solution is unstable, besides the two trivial eigenvalues and eigenfunctions, we need to discover at least one more eigenfunction with eigenvalue $\mu$ such that $\Re({\mu}) > 0$. By solving the linearized equation \eqref{eq:linear:op} via PINN for the second smooth solution, we do find another eigenfunction as shown in Figure \ref{fig:ccf_eig2}($c$) with associated eigenvalue $\mu = 0.36525$, which confirms that the second smooth solution at $\lambda = 0.6057$ is indeed an unstable solution. However, whether there exist more unstable eigenvalues for the second smooth solution requires systematic spectrum analysis, which is beyond the scope of the current paper and is our ongoing work.

As highlighted in the main text, such unstable solution has not been found by other conventional methods. The discovery of the unstable solution, as well as the universality of the PINN setting for various problems we investigated in this study demonstrates the powerfulness of PINN in discovering the self-similar blow-up solutions.

\section{Computational costs}
In this work, we use TensorFlow version 2.4.1 to perform the computation of the self-similar solution discovery via PINNs. All calculation is conducted under 64-bit double precision in order to guarantee a suitably small round-off error. It takes one million iterations to obtain the self-similar solution with the relative error of $10^{-4}$ for the 2D Bousinessq equation, which take 16 hours to compute on one NVIDIA A100 40GB GPU. For the 1D Burgers equation, it takes $100,000$ iterations to obtain the self-similar solution with the relative error of $10^{-7}$, which take less than 10 minutes to compute on one NVIDIA A100 40GB GPU. For the generalized De Gregorio equation, It takes one million iterations to obtain the self-similar solution with the relative error of $10^{-6}$, which takes 15 hours to compute on one NVIDIA A100 40GB GPU. The extra time spent on solving the generalized De Gregorio equation is due to the computation of high-accuracy numerical Hilbert Transform. 
\color{black}

\bibliographystyle{abbrv}
\bibliography{boussinesqPRL}